\documentclass[11pt]{amsart}
\usepackage{preamble}

\title[LCM of a Bivariate Quadratic Sequence]{The Least Common Multiple of a Bivariate Quadratic Sequence}

\author{Noam Kimmel}
\thanks{ This research was supported by the European Research Council (ERC) under the European Union's  Horizon 2020 research and innovation program  (Grant agreement No.    786758).}

\begin{document}

\maketitle 
\begin{abstract}

   Let $F\in\mathbb{Z}[x,y]$ be some polynomial of degree 2.
   In this paper we find the asymptotic behaviour of the least common multiple of the values of $F$ up to $N$.
   More precisely, we consider $\psi_F(N) = \log\left(\text{LCM}_{0<F(x,y)\leq N}\left\lbrace F(x,y)\right\rbrace\right)$ as $N$ tends to infinity.
   It turns out that there are 4 different possible asymptotic behaviours depending on $F$.
   For a generic $F$, we show that the function $\psi_F(N)$ has order of magnitude $\frac{N\log\log N}{\sqrt{\log N}}$.
   We also show that this is the expected order of magnitude according to a suitable random model.
   However, special polynomials $F$ can have different behaviours, which sometimes deviate from the random model.
   We give a complete description of the order of magnitude of these possible behaviours, and when each one occurs.

\end{abstract}

\setcounter{tocdepth}{1}
\tableofcontents
\newpage

\section{Introduction}
\subsection{The Prime Number Theorem}
Consider the Chebyshev function
$$
\psi(N) = \sum_{n\leq N}\Lambda(n)
$$
where $\Lambda$ is the von-Mangoldt function
\begin{equation*}
    \Lambda(n) = 
    \begin{cases}
    \log p& n=p^k \text{ for $p$ prime and }k>0\\
    0& \text{otherwise}
    \end{cases}.
\end{equation*}
The prime number theorem is equivalent to the assertion that $\psi(N)\sim N$.
An alternative definition for the function $\psi(N)$ can be given using the least common multiple of the first $N$ integers:
$$
\psi(N) = \log\left(\text{LCM}\left\lbrace1,2,3,4,....,N\right\rbrace\right).
$$
Thus, we see that the PNT can be stated in the following way:
\begin{theorem}[PNT]\label{PNT}
$$
\log\left(\mathop{\text{\normalfont  LCM}}_{n\leq N}\left\lbrace n\right\rbrace\right) \sim N.
$$
\end{theorem}
\subsection{LCM of consecutive polynomial values}
One can naturally wonder about generalizations of the PNT arising from the formulation given in \autoref{PNT}.
Let $f\in\mathbb{Z}[t]$ be some integer polynomial, and define the function
$$
\psi_f(N) = \log\left(
\mathop{\text{LCM}}_{n\leq N}
\left\lbrace f(n)\right\rbrace
\right).
$$
It is an interesting problem to try and understand the asymptotic behaviour of $\psi_f(N)$ for various $f$.

The case where $f$ is a linear polynomial was resolved in \cite{Bateman-1}.
This was later extended to all polynomials $f$ which are a product of linear terms in \cite{Hong-1}.
In both cases, it was shown that $\psi_f(N)\sim c_f N$ for some constant $c_f$.
The main ingredient of both proofs is Dirichlet's theorem on primes in arithmetic progressions.

In \cite{Cilleruelo-1}, Cilleruelo considers the case where $f$ is an irreducible polynomial of degree 2.
In this case he shows that 
$$
\psi_f(N)\sim N\log N.
$$
He further conjectures that if $f$ is an irreducible polynomial of degree $d>2$, then
$\psi_f(N)\sim(d-1)N\log N$.

While the upper bound 
$$
\psi_f(N)\lesssim (d-1)N\log N
$$
follows implicitly from Cilleruelo's work, proving the lower bound seems to be a difficult problem.
In \cite{Maynard-1} Maynard and Rudnick prove the lower bound 
$$
\psi_f(N)\gtrsim \frac{1}{d}N\log N,
$$
which is of the correct order of magnitude.
This was later improved by Sah in \cite{Sah-1} to give
$$
\psi_f(N)\gtrsim N\log N.
$$
However, while these bounds give the correct order of magnitude, there is still no irreducible polynomial of degree $d>2$ for which the conjecture is known to hold.
For further details regarding this problem see the survey \cite{Bazzanella-1}.

\subsection{LCM of multivariate polynomial values}
We consider another natural generalization of \autoref{PNT}.
Let $F\in \mathbb{Z}[x,y]$, and denote
$$
\psi_F(N) = \log\left(
\mathop{\text{LCM}}_{0<F(x,y)\leq N}
\left\lbrace F(x,y)\right\rbrace
\right).
$$
One can ask for the asymptotic behaviour of $\psi_F(N)$.
A similar problem was proposed by Shparlinski in \cite{Candela-1}.

For example, consider the case where $F(x,y)=x^2 + y^2$.
In this case, it is a classical result of Fermat that a number $n$ is of the form $x^2+y^2$ if and only if $n$ has the form
$$
n = 2^a\prod_{p_i\equiv 3 (\bmod 4)}p_i^{2b_i}
\prod_{p_j\equiv 1 (\bmod 4)}p_j^{c_j}.
$$
From which we get
$$
\psi_F(N) =
\left\lfloor \log_2 N\right\rfloor \log 2
+ 2\sum_{\substack{p^k\leq \sqrt{N}\\ p\equiv 3 (\bmod 4)}}{\log p}
+
\sum_{\substack{p^k\leq N\\ p\equiv 1 (\bmod 4)}}\log p.
$$
And so, from the prime number theorem in arithmetic progressions we get that $\psi_F(N)\sim\frac{1}{2}N$.
Note that this result relies heavily on the multiplicative properties of the image of $F$.
If we were to consider $F(x,y) = x^2+y^2+1$ instead, the method above wouldn't work.

\subsection{A Random model}
In order to better understand what to expect from $\psi_F(N)$ for $F(x,y) = x^2 + y^2 + 1$, we can consider the following random model.
Let $F\in\mathbb{Z}[x,y]$ be some polynomial, and denote
$$
\delta_F(N) = \frac{\#\left\lbrace
0<n\leq N \;\middle|\; n=F(x,y)
\right\rbrace}{N}.
$$
It was shown in \cite{Cilleruelo-2}, that if we take a random set $A\subset\lbrace 1,2,3,...,N \rbrace$ where each element is taken with probability $\delta = \delta(N)$, then 
$$
\psi(A) := \log\left(\mathop{\text{LCM}}_{a\in A}\lbrace a\rbrace \right) \sim N\frac{\log\left(\delta^{-1}\right) \delta}{1-\delta}
$$
almost surely as $N$ tends to infinity (provided $N \delta(N) \rightarrow \infty$).

And so, if we were to model the image of $F$ as a random subset of $\lbrace 1,2,3,4..,N\rbrace$ where each element is taken with probability $\delta_F(N)$, then we would expect 
\begin{equation*}
\psi_F(N)\sim N\frac{\log\left(\delta_F^{-1}(N)\right) \delta_F(N)}{1-\delta_F(N)}.
\end{equation*}
Assuming some local phenomena arising from small primes, we still expect to have
\begin{equation}\label{eq-i1}
\psi_F(N)\asymp N\frac{\log\left(\delta_F^{-1}(N)\right) \delta_F(N)}{1-\delta_F(N)}.
\end{equation}
We already know that this model doesn't predict the correct asymptotic for all $F$.
To see this, consider once more the example $F(x,y) = x^2 + y^2$.
In this case we have that $\delta_F(N) \asymp \frac{1}{\sqrt{\log N}}$, so the random model \eqref{eq-i1} would predict $\psi_F(N) \asymp \frac{N\log\log N}{\sqrt{\log N}}$.
But we saw that in this case $\psi_F(N) \sim \frac{1}{2}N$, which is significantly larger.

On the other hand, if we were to consider $F(x,y) = x^2 + y^2 + 1$ instead, then it turns out that $\psi_F(N)$ has order of magnitude $\frac{N \log \log N}{\sqrt{\log N}}$ which agrees with the order of magnitude predicted by the random model \eqref{eq-i1}.
This will follow from our main result in \autoref{main-thm3}.

The discrepancy between $\psi_F(N)$ for $F(x,y) = x^2 + y^2$ and the random model \eqref{eq-i1} can be attributed to the multiplicative structure of the image of $F$, which causes abnormally high correlations between primes and the image of $F$. 

\subsection{Main results}
We give a complete description of the order of magnitude of $\psi_F(N)$ where $F(x,y)\in\mathbb{Z}[x,y]$ is a polynomial of degree 2, and when it agrees with the random model \eqref{eq-i1}.
Denote
$$
F(x,y) = ax^2+bxy+cy^2+ex+fy+g\in\mathbb{Z}[x,y]
$$
and assume that $F$ represents arbitrarily large integers.
We denote 
$$\Delta = b^2-4ac$$ 
the discriminant of the homogeneous part of $F$, which we refer to as the discriminant of $F$.
We also denote 
$$D = af^2 + ce^2 -bfe + \Delta g$$ 
the Hessian determinant of the associated quadratic form 
$$
\overline{F}(x,y,z) = z^2 F(x/z,y/z).
$$
We refer to $D$ as the large discriminant of $F$.

As before, we also denote 
$$
\delta_F(N) = \frac{\#\left\lbrace
0<n\leq N \;\middle|\; n=F(x,y)
\right\rbrace}{N}.
$$

Our main result are the following theorems:
\begin{theorem}\label{main-thm1}
Let $F$ be as above, and assume that $\frac{\partial F}{\partial x},\frac{\partial F}{\partial y}$ are linearly dependent.
In this case we have that
$$
\delta_F(N) \sim \frac{1}{\sqrt{(a,c) N}}.
$$
\begin{enumerate}
\item If $F$ is irreducible then
$$
\psi_F(N) \sim N \frac{\log\left(\delta_F^{-1}(N)\right) \delta_F(N)}{1-\delta_F(N)} \sim \frac{\sqrt{N} \log N}{2 \sqrt{(a,c)}}.
$$
\item If $F$ is reducible then
$$
\psi_F(N) \sim c_F \sqrt{N}
$$
for some explicitly given constant $c_F$.
\end{enumerate}
\end{theorem}
\begin{remark}
The condition that $\frac{\partial F}{\partial x},\frac{\partial F}{\partial y}$ are linearly dependent means that the function $F$ is a function of essentially one variable (this claim will be made precise in the proof).
In this case, the random model \eqref{eq-i1} predicts the correct results when $F$ is irreducible, and even gives the correct constant, but fails in the reducible case.
\end{remark}

\begin{theorem}\label{main-thm2}
Let $F$ be as above, and assume that $\frac{\partial F}{\partial x},\frac{\partial F}{\partial y}$ are linearly independent.
\begin{enumerate}
\item If $\Delta$ is a perfect square, then $1 \ll \delta_F(N) \ll 1$ and 
$$
N \ll \psi_F(N) \ll N
$$
\item If $\Delta$ is not a perfect square and $D = 0$ then 
$\frac{1}{\sqrt{\log N}} \ll \delta_F(N) \ll \frac{1}{\sqrt{\log N}}$ and 
$$
N \ll \psi_F(N) \ll N.
$$
\end{enumerate}
\end{theorem}
\begin{remark}
The condition that $\Delta$ is a perfect square means that the homogeneous part of $F$ factors into linear terms over the rationals, and in this case the random model \eqref{eq-i1} predicts the correct order of magnitude.
The condition that $D=0$ means that $F$ is equivalent (over the rationals) to a (homogeneous) binary quadratic form. In this case the random model fails to give the correct order of magnitude for $\psi_F(N)$.
\end{remark}

\begin{theorem}\label{main-thm3}
Let $F$ be as above, and assume that $\frac{\partial F}{\partial x},\frac{\partial F}{\partial y}$ are linearly independent.
Assume also that $\Delta$ is not a perfect square, and that $D\neq 0$.
In this case $\frac{1}{\sqrt{\log N}} \ll \delta_F(N) \ll \frac{1}{\sqrt{\log N}}$, and
$$
\frac{N \log \log N}{\sqrt{\log N}} \ll \psi_F(N) \ll
\frac{N \log \log N}{\sqrt{\log N}}.
$$
\end{theorem}
\begin{remark}
In this case we have that
$$
\frac{N \log \log N}{\sqrt{\log N}} \asymp
\frac{\log\left(\delta_F^{-1}(N)\right) \delta_F(N)}{1-\delta_F(N)}
$$
which means that $\psi_F(N)$ has the same order of magnitude that we would expect from a random set according to \eqref{eq-i1}.

However, precise asymptotics seem to be beyond our method.
In \autoref{sec-6} we give a conjecture for the asymptotics of $\psi_{x^2 + y^2 + 1}(N)$, and present some numerical evidence.
\end{remark}

We handle the degenerate cases (\autoref{main-thm1}, \autoref{main-thm2}) in \autoref{sec-2}.
The bulk of the paper will be devoted to treating the generic case (\autoref{main-thm3}).
The main ingredients of the proof are the half dimensional sieve, the Bombieri-Vinogradov theorem, and genus theory for quadratic forms.

\subsection{Notations and Assumptions}
We say that a function $P(x,y)\in\mathbb{Z}[x,y]$ represents an integer $n$ if there exists $x_0,y_0\in\mathbb{Z}$ such that $n = P(x_0,y_0)$.
If there exists relatively prime $x_1,y_1\in\mathbb{Z}$ such that $n = P(x_1,y_1)$ we say that $n$ is properly represented by $P$.

Sums of the form $\sum_{n = F(x,y)\leq N}$ are to be regarded as summing over positive integers $n\leq N$ which are represented by $F$, without counting the number of different representations of $n$ as $F(x,y)$.

Unless otherwise stated, we use $(\cdot,\cdot)$ brackets to denote the GCD, and $[\cdot,\cdot]$ brackets to denote the LCM.
$F$ will denote the polynomial $ax^2 + bxy + cy^2 + ex + fy + g \in \mathbb{Z}[x,y]$.
We will denote $G_F(x,y) = ax^2 + bxy + cy^2$ the homogeneous part of $F$.
$\Delta = \Delta_F$ will denote the discriminant of $F$, and $D = D_F$ will denote the large discriminant of $F$.

When relevant, we assume all implied constants coming from the $\ll$ and $\mathcal{O}\left(\cdot\right)$ notations depend on $F$.

Throughout the paper, the letters $p,q$ will denote primes.

\section{Proof of the degenerate cases}\label{sec-2}
In this section we prove \autoref{main-thm1} and \autoref{main-thm2}.

\subsection{The single variable case}
We begin with a proof of \autoref{main-thm1}.
\begin{proof}
Let $F$ be as in \autoref{main-thm1}.
Since $\frac{\partial F}{\partial x},\frac{\partial F}{\partial y}$ are linearly dependent, there exists $A,B, \alpha,\beta,\gamma$ with $(A,B) = 1$ such that
$$
A \frac{\partial F}{\partial x} = B \frac{\partial F}{\partial y} = \alpha x + \beta y + \gamma.
$$
It follows that $\alpha A = \beta B$ and that
$$
F(x,y) = \frac{1}{2A\alpha}\left(\alpha x + \beta y + \gamma\right)^2 + C
$$
for some constant $C$.
And so, there is some quadratic polynomial $f(t) \in \mathbb{Z}[t]$ such that 
$$
F(x,y) = f\left(\frac{\alpha}{(\alpha,\beta)}x + \frac{\beta}{(\alpha,\beta)}y\right),
$$
with leading coefficient $(\alpha,\beta)^2/2A\alpha$.

Furthermore, we have that 
$$
a = \frac{\alpha}{2A} = \frac{\alpha^2}{2A\alpha},\quad  c = \frac{\beta}{2B} = \frac{\beta^2}{2A\alpha}.
$$
We then get that
$$
(a,c) = \left(\frac{\alpha^2}{2A\alpha}, \frac{\beta^2}{2A\alpha}\right) = \frac{(\alpha,\beta)^2}{2A\alpha}
$$
which is the leading coefficient of $f$.

To summarize, it can be seen that the image of $F(x,y)$ is equal to the image of $f(t)$, a quadratic polynomial with leading coefficient $(a,c)$.
From this we see that 
$$
\delta_F(N) \sim \frac{1}{\sqrt{(a,c)N}}.
$$

If $F$ is reducible, then so is $f$.
This means that $f$ is a product of linear terms.
Thus, from \cite{Hong-1} we get that
$$
\psi_F(N) = 
\log\left(\mathop{\text{LCM}}_{0<n \leq \sqrt{N/(c,a)}} 
\left\lbrace f(n) \right\rbrace \right) \sim c_f \frac{\sqrt{N}}{\sqrt{(c,a)}}
$$
for some explicitly given constant $c_f$.

Assume now that $F$ is irreducible.
This means that $f$ is also irreducible.
Thus, from \cite{Cilleruelo-1} we get that
$$
\psi_F(N) = 
\log\left(\mathop{\text{LCM}}_{0<n \leq \sqrt{N/(c,a)}} 
\left\lbrace f(n) \right\rbrace \right) \sim 
\frac{\sqrt{N}\log N}{2 \sqrt{(a,c)}}.
$$
\end{proof}

\subsection{The degenerate discriminants case}
Before proving \autoref{main-thm2} we prove the following two lemmas.
\begin{lemma}\label{lem-1}
Let $F$ be as in \autoref{main-thm2}.
If the image of $F$ contains an arithmetic progression, then 
$$
N\ll \psi_F(N) \ll N.
$$
\end{lemma}
\begin{proof}
The upper bound $\psi_F(N) \ll N$ always holds since
$$
\psi_F(N) \ll 
\log\left(\text{LCM}\left\lbrace 1,2,3,4,..., N\right\rbrace\right) \sim N.
$$
Assume now that $\text{Im}(F)$ contains an arithmetic progression $\left\lbrace An + B\right\rbrace_{n\in\mathbb{Z}}$, $A\neq 0$.
Denote $G = (A,B)$, and $A' = A/G,\, B' = B/G$.
Note that $(A',B')=1$.
In this case we have
\begin{multline*}
\psi_F(N) \gg    
\log\left(\mathop{\text{LCM}}_{0<n\leq (N-B)/A}
\left\lbrace A' n + B'\right\rbrace \right)
\gg \\
\sum_{0<n\leq (N-B)/A}\Lambda\left(A' n + B'\right) \gg N
\end{multline*}
where the last inequality follows from Dirichlet's theorem on primes in arithmetic progressions. 
\end{proof}

We shall say that a polynomial $P(x,y) \in \mathbb{Z}[x,y]$ belongs to the class $\mathcal{H}$ if for every integer $A\neq 0$, $P$ represents an integer prime to $A$.
\begin{lemma}\label{lem-2}
Let 
$$
F = ax^2 + bxy + cy^2 + ex + fy + g\in\mathbb{Z}[x,y].
$$
Assume that $F$ represents arbitrarily large integers, and that $\partial F / \partial x, \partial F / \partial y$ are linearly independent.
There exists $H\in \mathcal{H}$ and a constant $C_F$, such that $H$ has linearly independent partial derivatives, represents arbitrarily large integers, and the following holds:
$$
\delta_H(N) \ll \delta_F(C_F N), \quad 
\psi_H(N) \ll \psi_F(C_F N).
$$
Furthermore, the small and large discriminants of $H$ satisfy
$$
\Delta_H = C_F^{-2}\Delta_F, \quad D_H =  C_F^{-2}D_F.
$$
\end{lemma}

\begin{proof}
Denote 
$$
F(x,y) = ax^2 + bxy + cy^2 + ex + fy + g.
$$
From \cite[Lemma 2]{Iwaniec-1} we get that the condition that $F$ be in the class $\mathcal{H}$ is equivalent to the conditions
\begin{gather*}
    (a,b,c,e,f,g)=1 \\
    a,c,b,g \not\equiv e,f,0,0\; (\bmod 2).
\end{gather*}

Denote $W = (a,b,c,e,f,g)$.
We start by defining $G = W^{-1}F$.
Denote 
$$
G(x,y) = a'x^2 + b'xy + c'y^2 + e'x + f'y + g'.
$$
It is easy to see that $G$ still represents arbitrarily large integers, and that 
$$
\psi_{F}(W N) \sim \psi_{G}(N),\quad
\delta_{F}(W N) \sim \delta_{G}(N)
$$
If $G\in\mathcal{H}$, pick $H=G$ and $C_F = W$.
Otherwise, we must have that $G$ represents only even integers.
We then define
$H(x,y) = \frac{1}{2}G(2x,2y)$.
From \cite[Lemma 2]{Iwaniec-1} we see that $H$ is in $\mathcal{H}$.
Furthermore, it can be easily verified that $H$ has linearly independent partial derivatives, and represents arbitrarily large integers.
It can also be seen that
$$
\psi_{H}(N) \ll \psi_G(2N),\quad 
\delta_{H}(N) \ll \delta_G(2N)
$$
And so, the statement of the lemma holds with $H$ and $C_F = 2W$.
The statement about the discriminants of $H$ follow from computation.
\end{proof}

We now give a proof of \autoref{main-thm2}.
\begin{proof}
Let $F$ be as in  \autoref{main-thm2}, and assume that $\partial F/\partial x$, $\partial F/\partial y$ are linearly independent.

We begin by considering the case in which $\Delta\neq 0$ is a perfect square.
without loss of generality, we can assume $a\neq 0$.
In this case, the quadratic part of $F$, which we denote
$$
G_F(x,y) = ax^2 + bxy + cy^2,
$$
satisfies $4aG_F(x,y) = (2a x + r_1 y ) (2a x + r_2 y )$ with $r_1,r_2 = b \pm \sqrt{\Delta}$.
Since $\Delta\neq 0$, at least one of the vectors $(2a,r_1),(2a,r_2)\in\mathbb{Z}^2$ is not proportional to $(f,e)\in\mathbb{Z}^2$.
Assume without loss of generality that $(2a,r_1)$ satisfies this requirement.
We then have 
$$
F(r_1t, -2at)
=
t\cdot  (r_1e -2af) + g.
$$
It follows that the image of $F$ contains an arithmetic progression.
Thus, from \autoref{lem-1} we get that 
$$
N\ll \psi_F(N) \ll N.
$$

Assume now that $\Delta = 0$.
In this case the quadratic part of $F$ is a square.
Assuming once more that $a\neq 0$ we get that
$$
4aF(x,y) 
= (2ax + by)^2 + 2e(2ax + by) + 2(2af - eb) y + 4ag.
$$
Since $\frac{\partial F}{\partial x}, \frac{\partial F}{\partial y}$ are linearly independent, we must have that $2af - eb\neq 0$.
It follows that 
$$
F(bt,-2at) = t\cdot (2af - be) + g.
$$
Thus, the image of $F$ contains an arithmetic progression, so from \autoref{lem-1} we get 
$$
N\ll \psi_F(N) \ll N.
$$

Note that we showed that if $\Delta$ is a perfect square (zero or non-zero), the image of $F$ contains an arithmetic progression.
And so, in this case we have
$$
1 \ll \delta_F(N)  \ll 1.
$$

Consider now the case where $\Delta$ is not a perfect square, and $D = 0$.
Let $H\in\mathbb{Z}[x,y]$ be as in \autoref{lem-2}, and $C_F$ the constant from the lemma.
Like $F$, the discriminant of $H$ is not a perfect square, and the large discriminant is 0.
From \cite[Proposition 3]{Iwaniec-1}, we get that
$$
\sum_{\substack{p\leq N \\ p \text{ prime} \\ p = H(x,y)}} 1 
\gg \frac{N}{\log N}.
$$
It follows that $\psi_{H}(N) \gg N$.
Since
$$
\psi_{F}(N) \gg \psi_{H}\left(C_F N\right) \gg N
$$ 
we get that $\psi_{F}(N) \gg N$ as required.
Note that the upper bound $\psi_{F}(N) \ll N$ is always true as previously noted.

The bounds
$$
\frac{1}{\sqrt{\log N}}\ll \delta_F(N) \ll \frac{1}{\sqrt{\log N}}
$$
are given in \autoref{lem-47} and \autoref{lem-48}.
\end{proof}

\section{Proof of the generic case}
In this section we prove \autoref{main-thm3}.
That is, we show that for $F$ satisfying the assumption of \autoref{main-thm3},
$$
\frac{N\log\log N }{\sqrt{\log N}} \ll \psi_F(N) \ll \frac{N\log\log N }{\sqrt{\log N}}.
$$
For the proof we will need estimations for the number of small multiples of large primes represented by $F$.
These estimates will be stated in \autoref{lem-43}, \autoref{lem-44}, \autoref{lem-45}, which we will prove in \autoref{sec-4}.
The claim that $\delta_F(N)$ has order of magnitude $\frac{1}{\sqrt{\log N}}$ is given in \autoref{lem-47} and \autoref{lem-48}.

Throughout the section we assume that 
$$
F(x,y) = ax^2 + bxy + cy^2 + ex + fy + g\in\mathbb{Z}[x,y]
$$
satisfies the assumptions of \autoref{main-thm3}.
We denote $\Delta = b^2 -4ac$ the discriminant of $F$, and $D = af^2 + ce^2 -bfe + \Delta g$ the large discriminant of $F$.

\subsection{The Upper Bound}
In this section we prove the following:
\begin{proposition}\label{prop-1}
$$
\psi_F(N) \ll \frac{N\log\log N}{\sqrt{\log N}}.
$$
\end{proposition}
We begin by stating some preliminary lemmas, and then show how
\autoref{prop-1} follows.
\begin{lemma}\label{lem-41}
Denote
$$
\mathop{\text{\normalfont LCM}}_{0<F(x,y)\leq N}
\left\lbrace F(x,y)\right\rbrace
 = \prod_{\substack{p\leq N \\ p \text{ prime}}}p^{\theta_p}.
$$
Then
$$
\sum_{\substack{p\leq N \\ p \text{ prime} \\ \theta_p \geq 2}}
(\theta_p - 1)\log p = \mathcal{O}\left(\sqrt{N}\right).
$$
\end{lemma}
This lemma shows that the contribution to $\psi_F(N)$ coming from multiplicities is negligible.
\begin{proof}
$$
\sum_{\substack{p\leq N \\ p \text{ prime} \\ \theta_p > 1}}
(\theta_p - 1)\log p
\leq
\sum_{m\geq 2}\;
\sum_{\substack{p\leq N^{1/m} \\ p \text{ prime}}} \log p \ll \sqrt{N}.
$$
\end{proof}
\begin{lemma}\label{lem-42}
With the notations from \autoref{lem-41}, and for $s>0$
$$
\sum_{\substack{p\leq N/\log^{s} N \\ p \text{ prime} \\ \theta_p \geq 1}}\log p \ll \frac{N}{\log^s N}.
$$
\end{lemma}
This lemma will be used to bound the contribution of small primes to the LCM.
\begin{proof}
$$
\sum_{\substack{p\leq N/\log^{s} N \\ p \text{ prime} \\ \theta_p \geq 1}}\log p \leq 
\sum_{\substack{p\leq N/\log^{s} N \\ p \text{ prime}}}\log p \ll \frac{N}{\log^s N}.
$$
\end{proof}

\begin{lemma}\label{lem-43}
For $k\leq \log^5 N$ denote
$$
S_k(N) = \sum_{\substack{  kp = F(x,y) \leq N \\p \text{ prime}}} 1.
$$
Then
$$
S_k(N) \ll \frac{N}{k \log^{3/2} N}
\prod_{\substack{q \mid k \\ \left(\frac{\Delta}{q}\right) = -1}}\left(1 - \frac{1}{q}\right)^{-1}.
$$
\end{lemma}
We postpone the proof of this lemma to \autoref{sec-4}.

We can now prove the upper bound for $\psi_F(N)$.
\begin{proof}[Proof of \autoref{prop-1}]
from \autoref{lem-41}, \autoref{lem-42} we get that
\begin{equation}\label{eq-41}
\psi_F(N) = \sum_{\substack{ p \text{ prime} \\ p \mid F(x,y) \leq N\\ p \geq N \log^{-1/2 }N }} \log p
 \;+\; \mathcal{O}\left(\frac{N}{\sqrt{\log N}}\right).
\end{equation}
With the notation of \autoref{lem-43} we get that
$$
\sum_{\substack{ p \text{ prime} \\ p \mid F(x,y) \leq N\\ p \geq N \log^{-1/2 }N }} \log p \; \ll \;
\log N \sum_{k\leq \sqrt{\log N}}S_k(N).
$$
From \autoref{lem-43} we then get
\begin{multline}\label{eq-42}
\sum_{\substack{ p \text{ prime} \\ p \mid F(x,y) \leq N\\ p \geq N \log^{-1/2 }N }} \log p \; \ll \; \\  
\log N
\frac{N}{\log^{3/2} N}
\sum_{k \leq \sqrt{\log N}}
\frac{1}{k}\prod_{\substack{q \mid k \\ \left(\frac{\Delta}{q}\right) = -1}}\left(1 - \frac{1}{q}\right)^{-1}
\;\ll\;
\frac{N \log \log N}{\sqrt{\log N}}.
\end{multline}
And so, from \eqref{eq-41},\eqref{eq-42} we get that
$$
\psi_F(N) \ll \frac{N \log \log N}{\sqrt{\log N}}.
$$
\end{proof}

\subsection{The Lower Bound}
We now look at the lower bound.
\begin{proposition}\label{prop-2}
$$
\psi_F(N) \gg \frac{N\log\log N}{\sqrt{\log N}}.
$$
\end{proposition}

We begin by stating some relevant lemmas.
\begin{lemma}\label{lem-44}
For $k\leq \log^5 N$ denote
$$
S_k(N) = \sum_{\substack{  kp = F(x,y) \leq N \\p \text{ prime}}} 1.
$$
There exists a constant $c_F$ depending only on $F$, such that if $k$ is prime to $c_F$ then
$$
S_k(N) \gg \frac{N}{k \log^{3/2} N}
\prod_{\substack{q \mid k \\ \left(\frac{\Delta}{q}\right) = -1}}\left(1 - \frac{1}{q}\right)^{-1}
$$
\end{lemma}
\begin{lemma}\label{lem-45}
Let $k_2 < k_1 \leq \log^5 N$, and denote
$$
S_{k_1,k_2}(N) = \sum_{\substack{ p \text{ prime} \\ k_1p = F(x_1,y_1) \leq N \\ k_2p = F(x_2,y_2) \leq N}} 1.
$$
Then
$$
S_{k_1,k_2}(N) \ll \frac{N}{k_1 \log^{2} N}
\prod_{\substack{q \mid k_1 k_2 (k_1 - k_2) \\ q>3}}\left(1 - \frac{3}{q}\right)^{-1}
$$
\end{lemma}
We postpone the proofs of \autoref{lem-44} and \autoref{lem-45} to \autoref{sec-4}.

We now prove the lower bound for $\psi_F(N)$:
\begin{proof}[proof of \autoref{prop-2}]
By ignoring small primes, we have that
$$
\psi_F(N) \geq 
\sum_{\substack{ p \text{ prime} \\ p \mid F(x,y) \leq N\\ p \geq N \log^{-1/2 +\epsilon}N }} \log p
\;\sim\;
\log N
\sum_{\substack{ p \text{ prime} \\ p \mid F(x,y) \leq N\\ p \geq N \log^{-1/2 +\epsilon}N }} 1.
$$
An iteration of the inclusion-exclusion principle then gives
\begin{multline}\label{eq-43}
\frac{1}{\log N}\psi_F(N) \gg \\ 
\sum_{\substack{k \leq \log^{1/2 - \epsilon} N \\ (k,c_F)=1}}
S_k(N) \;
 - 
\sum_{ k_1 \leq \log^{1/2 - \epsilon} N}\;
\sum_{ k_2  < k_1} S_{k_1,k_2}(N)
\end{multline}
where $c_F$ is the constant from \autoref{lem-44}.

For the first sum we use \autoref{lem-44} and we get
\begin{equation}\label{eq-44}
\sum_{\substack{k \leq \log^{1/2 - \epsilon} N \\ (k,c_F)=1}}
S_k(N)
\gg \frac{N}{\log^{3/2} N} \sum_{\substack{k\leq \log^{1/2 - \epsilon} N \\ (k,c_F)=1}} \frac{1}{k}
\gg \frac{N \log \log N}{\log^{3/2} N}.
\end{equation}

For the second double sum we use \autoref{lem-45}, noting that 
$$\prod_{\substack{q \mid k_1 k_2 (k_1 - k_2) \\ q>3}}\left(1 - \frac{3}{q}\right)^{-1} \ll (\log \log N)^3,
$$
we get
\begin{multline}\label{eq-45}
\sum_{ k_1 \leq \log^{1/2 - \epsilon} N}\;
\sum_{ k_2  < k_1} S_{k_1,k_2}(N) \ll \\
\sum_{k_1 \leq \log^{1/2 - \epsilon} N} \frac{1}{k_1}
\sum_{ k_2  < k_1} \frac{N (\log\log N)^3}{\log^2 N }
\ll \frac{N (\log \log N)^3}{\log^{3/2 + \epsilon} N}.
\end{multline}

From \eqref{eq-43}, \eqref{eq-44}, \eqref{eq-45} we get 
$$
\psi_F(N) \gg \frac{N \log\log N}{\sqrt{\log N}}
$$
as required.
\end{proof}


\section{preliminary Lemmas}
In this section we state and prove some lemmas which will be used in the following section.

\subsection{Lemmas From Sieve Theory}
\begin{theorem}[Bombieri-Vinogradov Theorem]\label{b-v}
$$
\sum_{Q\leq \frac{\sqrt{N}}{\log^{35}N} }\;
\max_{\substack{a\\ (a,Q) = 1}}\left| \pi(N,Q;a) - \frac{\text{Li} (N)}{\phi(Q)}\right| \ll \frac{N}{\log^{30} N}.
$$
\end{theorem}

\begin{proof}
Follows from \cite[Theorem 1]{Granville}
\end{proof}

\begin{theorem}[Main theorem of the $\beta$-sieve]\label{thm-sieves}
Let $\kappa \geq 0$ be a real number (referred to as the sieve dimension).
Let $\mathcal{P}$ be a set of primes.
Let $\mathcal{A} = (a_n)$ be a finite sequence of non-negative numbers.
Denote $\mathcal{A}_\lambda$ the subsequence of $\mathcal{A}$ consisting of those $a_n$ with $\lambda \mid n$.
We denote also $|\mathcal{A}_\lambda | = \sum_{\lambda \mid n} a_n$.

Let $X>0$ be a real number, and let $g(\lambda)$ be a multiplicative function satisfying $0\leq g(p) < 1$.
Denote
$$
\left|\mathcal{A}_\lambda\right| = g(\lambda) X + r_\lambda.
$$
Assume that there exists $L\geq 1$ such that for all $2\leq w < z$ we have
\begin{equation}\label{eq:g-req}
\prod_{\substack{w\leq p < z \\ p \in \mathcal{P}}}\left(1 - g(p)\right)^{-1}
\leq \left(\frac{\log z}{\log w}\right)^{\kappa}\left(1 + \frac{L}{\log w}\right).
\end{equation}
Let $2<z\leq D$, and denote $s = \log D / \log z$.
Denote $P(z) = \prod_{\substack{p< z \\ p\in\mathcal{P}}} p$ and $\mathcal{S}(\mathcal{A},z) = \sum_{(n,P(z)) = 1}a_n$.
Then we have that
\begin{equation}\label{eq:sieve1}
\mathcal{S}(\mathcal{A},z) \leq 
X \prod_{\substack{p<z \\ p\in\mathcal{P}}}(1-g(p))\left\lbrace C_\kappa(s) + \mathcal{O}\left(\frac{1}{\log^{1/6}D}\right)\right\rbrace + 
\sum_{\substack{\lambda \mid P(z) \\ \lambda \leq D}}|r_\lambda|,
\end{equation}
\begin{equation}\label{eq:sieve2}
\mathcal{S}(\mathcal{A},z) \geq 
X \prod_{\substack{p<z \\ p\in\mathcal{P}}}(1-g(p))\left\lbrace c_\kappa(s) + \mathcal{O}\left(\frac{1}{\log^{1/6}D}\right)\right\rbrace - 
\sum_{\substack{\lambda \mid P(z) \\ \lambda \leq D}}|r_\lambda|,
\end{equation}
where the functions $C_\kappa,c_\kappa$ satisfy the following.
\begin{equation}\label{eq:sieve3}
C_\kappa(s) = 1 + \mathcal{O}\left(e^{-s}\right),\quad
c_\kappa(s) = 1 + \mathcal{O}\left(e^{-s}\right).
\end{equation}
If $\kappa = 1/2$ then $s^{1/2}C(s) = 2\sqrt{e^\gamma / \pi}$ for $0\leq s \leq 2$, and
\begin{equation*}
s^{1/2}c(s) = \sqrt{\frac{e^\gamma}{\pi}}\int_{1}^{s}\frac{dt}{(t(t-1))^{1/2}}
\end{equation*}
for $1\leq s \leq 3$.

The constants in the $\mathcal{O}$ notations depend only on $\kappa,L$.
\end{theorem}
\begin{proof}
The bounds \eqref{eq:sieve1}, \eqref{eq:sieve2} are given in \cite[Theorem 11.13]{Friedlander-1}.
The claim in \eqref{eq:sieve3} is given in \cite[Equation 11.68]{Friedlander-1}.
The expressions for $C_{1/2},c_{1/2}$ are given in \cite[Chapter 14]{Friedlander-1}

\end{proof}

We are going to use \autoref{thm-sieves} with $\kappa = 1/2$ and $\kappa = 2$.
When applied with $\kappa = 2$, we will also need a slight variation of \autoref{thm-sieves} which allows us to sift by non-zero residue classes.

\begin{theorem}\label{thm:sieves-variant}
Let $\mathcal{P}$ be a set of primes.
Let $\mathcal{A} = (a_n)$ be a finite sequence of non-negative numbers.
For each $p\in \mathcal{P}$ let $\Omega_p \subsetneqq \ZZ/p\ZZ$ be a set of residue classes mod $p$.
Let $g$ be a multiplicative function with $0\leq g(p) < 1$.
Let $\kappa \geq 0$ be some number, and assume $g$ satisfies \eqref{eq:g-req} with $\kappa$ and some $L\geq 1$.
Denote $\mathcal{A}^\Omega_\lambda$ the subsequence of $\mathcal{A}$ consisting of those $n$'s for which $n\bmod p$ is in $\Omega_p$ for all $p\mid \lambda$, and denote 
$$
|\mathcal{A}^\Omega_\lambda| = \sum_{\substack{n\\\forall p\mid \lambda : \; n\bmod p\in\Omega_p} }a_n.
$$
Let $X>0$ be a real number, and denote 
$$
|\mathcal{A}^\Omega_\lambda| = Xg(\lambda) + r_\lambda.
$$

Let $2<z\leq D$, and denote $s = \log D / \log z$, $P(z) = \prod_{\substack{p< z \\ p\in\mathcal{P}}} p$.
We further denote 
$$
\mathcal{S}\left(\mathcal{A},z,\Omega\right) = \sum_{\substack{n \\ \forall p\mid P(z) : \; n\bmod p \not \in \Omega_p}} a_n.
$$
Then
\begin{equation*}
\mathcal{S}(\mathcal{A},z, \Omega) \leq 
X \prod_{\substack{p<z \\ p\in\mathcal{P}}}(1-g(p))\left\lbrace C_\kappa(s) + \mathcal{O}\left(\frac{1}{\log^{1/6}D}\right)\right\rbrace + 
\sum_{\substack{\lambda \mid P(z) \\ \lambda \leq D}}|r_\lambda|,
\end{equation*}
\begin{equation*}
\mathcal{S}(\mathcal{A},z, \Omega) \geq 
X \prod_{\substack{p<z \\ p\in\mathcal{P}}}(1-g(p))\left\lbrace c_\kappa(s) + \mathcal{O}\left(\frac{1}{\log^{1/6}D}\right)\right\rbrace - 
\sum_{\substack{\lambda \mid P(z) \\ \lambda \leq D}}|r_\lambda|,
\end{equation*}
where the functions $C_\kappa,c_\kappa$ satisfy
\begin{equation*}
C_\kappa(s) = 1 + \mathcal{O}\left(e^{-s}\right),\quad
c_\kappa(s) = 1 + \mathcal{O}\left(e^{-s}\right).
\end{equation*}
The constants in the $\mathcal{O}$ notations depend only on $\kappa,L$.
\end{theorem}
\begin{proof}
Define the function 
$$
f(n) = \prod_{\substack{p\in \mathcal{P}\\ p \leq z \\ n\bmod p \in \Omega_p}} p.$$
We now define the sequence $\mathcal{B} = (b_m)$ by
$$
b_m = \sum_{f(n) = m}a_n.
$$
With these notations, $|\mathcal{A}^\Omega_\lambda| = |\mathcal{B}_\lambda|$ and $\mathcal{S}(\mathcal{A},z,\Omega) = \mathcal{S}(\mathcal{B},z)$.
And so, applying \autoref{thm-sieves} to the sequence $\mathcal{B}$ gives the required result.
\end{proof}

\begin{remark}
We will at times apply \autoref{thm-sieves} and \autoref{thm:sieves-variant} with a finite set $\mathcal{M}\subset \NN$ in place of the sequence $\mathcal{A}$, in which case the sequence $\mathcal{A}$ is to be understood as the indicator series of $\mathcal{M}$: $a_n = \mathbb{1}_{\mathcal{M}}(n)$.
\end{remark}

\subsection{Technical Lemmas}
\begin{lemma}\label{lem-a3}
Let $A\neq 0$ be an integer.
If an integer $\Delta\equiv 0,1 (\bmod 4)$ is different from a perfect square then there exists a constant $c = c(\Delta)$ which does not depend on $A$ such that for all $x>2$:
\begin{multline*}
\sum_{\substack{m\leq x \\ q\mid m \Rightarrow \left(\frac{\Delta}{q}\right)=1 \\ (m,A)=1}}m^{-1}\log^{-2}\left(3 + \frac{x}{m}\right)\prod_{q\mid m}\left(1-\frac{1}{q}\right)^{-2}
\leq \\
\leq c \log^{-3/2}x \prod_{\substack{q\mid A \\ \left(\frac{\Delta}{q}\right) = 1}}\left(1 - \frac{1}{q}\right).
\end{multline*}
\end{lemma}
\begin{proof}
Throughout the proof, the constants $c_i$ will depend on $\Delta$.

We have that
\begin{multline*}
\sum_{\substack{m \leq x \\  q\mid m \Rightarrow \left(\frac{\Delta}{q}\right) = 1 \\ (m, A) = 1}} 
\frac{1}{m}
\prod_{q\mid m}\left(1 - \frac{1}{q}\right)^{-2}
\leq \\
\prod_{\substack{q\leq x \\ \left(\frac{\Delta}{q}\right) = 1 \\ (q,A) = 1}} \left(1 + \frac{1}{q} + \frac{2}{q(q-1)} + \frac{1}{q (q-1)^2}\right)
\leq
c_1 \log^{1/2}x \prod_{\substack{q\mid A \\ \left(\frac{\Delta}{q}\right) = 1}}\left(1 - \frac{1}{q}\right).
\end{multline*}
Denote $a(m) = \log^{-2}\left(3 + \frac{x}{m}\right)$.
Then for $1\leq m\leq x$:
$$
a'(m) \leq 2 m^{-1} \log^{-3}\left(3+\frac{x}{m}\right).
$$
By partial summation we now get
\begin{multline*}
\sum_{\substack{m\leq x \\ q\mid m \Rightarrow \left(\frac{\Delta}{q}\right)=1 \\ (m,A)=1}}
m^{-1}\log^{-2}\left(3 + \frac{x}{m}\right)\prod_{q\mid m}\left(1-\frac{1}{q}\right)^{-2}
\leq \\
c_2\prod_{\substack{q\mid A \\ \left(\frac{\Delta}{q}\right) = 1}}\left(1 - \frac{1}{q}\right) 
\int_{m=2}^{x} \frac{\log^{1/2}(m)}{m\log^{3}\left(3 + \frac{x}{m}\right)}dm
\leq 
c_3 \log^{-3/2}x \prod_{\substack{q\mid A \\ \left(\frac{\Delta}{q}\right) = 1}}\left(1 - \frac{1}{q}\right)
\end{multline*}
\end{proof}

\begin{lemma}\label{lem-a4}
Let $\Delta \neq 0$ be some integer.
Then for $X \geq 1$
\begin{equation}\label{eq-a1}
\sum_{\substack{r_1 \\ (r_1,\Delta) = 1}}
\sum_{\substack{r_2 \\ (r_2,\Delta) = 1}}
\sum_{\substack{d_1 \\ q\mid d_1\Rightarrow q\mid \Delta}}
\sum_{\substack{d_2 \\ q\mid d_2\Rightarrow q\mid \Delta \\ [r_1^2 d_1, r_2^2 d_2] \geq X}} \frac{1}{[r_1^2 d_1, r_2^2 d_2]}
\leq c(\Delta) \frac{(\log X)^2 + 1}{\sqrt{X}} 
\end{equation}
\end{lemma}
\begin{proof}
Throughout the proof, the constants $c_i$ will depend on $\Delta$.

We denote
$$
\tau(n) = \sum_{[a,b]=n}1 = \sum_{d\mid n^2} 1.
$$
For $s>0$ we have
\begin{equation}\label{eq-a2}
\sum_{k=0}^{\infty}\frac{\tau(p^k)}{p^{sk}}=
\sum_{k=0}^{\infty}(2k+1)p^{-sk}
= \frac{1 + p^{-s}}{(1 - p^{-s})^2}.
\end{equation}
It is also known that
\begin{equation}\label{eq-a3}
\sum_{n\leq N}\tau(n) \asymp N\log^2 N.
\end{equation}
By changing the order of summation, enumerating over $l_r = [r_1,r_2]$ and $l_d = [d_1,d_2]$, we get that the LHS of \eqref{eq-a1}, which we will denote $\Sigma$, satisfies
$$
\Sigma = 
\sum_{\substack{l_r \\ (l_r,\Delta) =1 }}\;
\sum_{\substack{l_d \\ q\mid l_d \Rightarrow q\mid \Delta  \\ l_r^2 l_d \geq X}}
\frac{\tau(l_r)\tau(l_d)}{l_r^2 l_d}.
$$
We separate the last sum into two sums as follows:
\begin{multline}\label{eq-a4}
\Sigma \leq \\
\sum_{\substack{l_d < X \\ q\mid l_d \Rightarrow q\mid \Delta }}\;
\sum_{\substack{l_r \geq \left(X / l_d\right)^{1/2} \\ (l_r,\Delta) =1 }}
\frac{\tau(l_r)\tau(l_d)}{l_r^2 l_d}
 \;+\;
\sum_{\substack{l_d \geq X \\ q\mid l_d \Rightarrow q\mid \Delta }}
\sum_{\substack{l_r \\ (l_r,\Delta) =1 }}
\frac{\tau(l_r)\tau(l_d)}{l_r^2 l_d}
\;=\; \Sigma_1 + \Sigma_2.
\end{multline}

We begin by estimating $\Sigma_1$.
From \eqref{eq-a3} and partial summation we have that
$$
\sum_{\substack{l_r > \left(X / l_d\right)^{1/2} \\ (l_r,\Delta) =1 }}
\frac{\tau(l_r)}{l_r^2}
\leq c_1 \frac{\log^2 X}{\sqrt{X}}\sqrt{l_d}.
$$
And so, we get that
\begin{multline*}
\Sigma_1 \leq 
\sum_{\substack{l_d < X \\ q\mid l_d \Rightarrow q\mid \Delta }}
\frac{\tau(l_d)}{l_d}
\sum_{\substack{l_r \geq \left(X / l_d\right)^{1/2} \\ (l_r,\Delta) =1 }}
\frac{\tau(l_r)}{l_r^2}
\leq c_2
\frac{\log^2 X}{\sqrt{X}}
\sum_{\substack{l_d\\ q\mid l_d \Rightarrow q\mid \Delta }}
\frac{\tau(l_d)}{\sqrt{l_d}}
=\\
c_2\frac{\log^2 X}{\sqrt{X}}
\prod_{q\mid\Delta}\frac{1 + q^{-1/2}}{(1 - q^{-1/2})^2}
\end{multline*}
where for the last equality we used \eqref{eq-a2}.
It follows that 
\begin{equation}\label{eq-a5}
\Sigma_1 \leq  c_3 \frac{\log^2 X}{\sqrt{X}}.    
\end{equation}

As for $\Sigma_2$ we have that for $\delta>0$:
\begin{multline*}
\Sigma_2 \leq c_4
\sum_{\substack{l_r}}\frac{\tau(l_r)}{l_r^2}
\sum_{\substack{l_d > X \\ q\mid l_d \Rightarrow q\mid \Delta }}
\frac{\tau(l_d)}{l_d}
\leq c_5 
\sum_{\substack{l_d \\ q\mid l_d \Rightarrow q\mid \Delta }}
\frac{\tau(l_d)}{l_d}\left(\frac{l_d}{X}\right)^\delta
= \\ 
\frac{c_5}{X^\delta}
\sum_{\substack{l_d \\ q\mid l_d \Rightarrow q\mid \Delta }}
\frac{\tau(l_d)}{l_d^{1 - \delta}}
= \frac{c_5}{X^\delta}
\prod_{q\mid \Delta}\left(\frac{1 + q^{-1 + \delta}}{(1 - q^{-1+\delta})^2}\right) 
\leq \frac{c_6(\delta)}{X^\delta}
\end{multline*}
where we once more used \eqref{eq-a2}.
Choosing $\delta = 1/2$ we get
\begin{equation}\label{eq-a6}
\Sigma_2 \ll X^{-1/2}.    
\end{equation}

Thus, from \eqref{eq-a3},\eqref{eq-a5}, \eqref{eq-a6} we get
$$
\Sigma \leq c\frac{(\log^2 X) + 1}{\sqrt{X}}
$$
as required.
\end{proof}

\begin{lemma}\label{lem-a5}
Let $\Delta \neq 0$ be some integer.
Then for $X > 0$
\begin{equation}\label{eq-a7}
\sum_{\substack{r_1 \\ (r_1,\Delta) = 1}}
\sum_{\substack{r_2 \\ (r_2,\Delta) = 1}}
\sum_{\substack{d_1 \\ q\mid d_1\Rightarrow q\mid \Delta}}
\sum_{\substack{d_2 \\ q\mid d_2\Rightarrow q\mid \Delta \\ [r_1^2 d_1, r_2^2 d_2] \leq X}} 1
\leq c(\Delta) \sqrt{X} \log^2 X
\end{equation}
\end{lemma}
\begin{proof}
We use the same notations from \autoref{lem-a4}.
Denoting the LHS of \eqref{eq-a7} by $\Sigma$, we get that
\begin{equation*}
\Sigma \leq 
\sum_{\substack{l_d \leq X\\ q\mid l_d \Rightarrow q\mid \Delta}}
\tau(l_d)
\sum_{\substack{l_r \leq \sqrt{X/l_d}}} \tau(l_r).
\end{equation*}
From \eqref{eq-a3} we then get that
$$
\Sigma \leq \sqrt{X}\log^2 X 
\sum_{\substack{l_d \\ q\mid l_d \Rightarrow q\mid \Delta}}
\frac{\tau(l_d)}{\sqrt{l_d}}
= \sqrt{X}\log^2 X \prod_{q\mid \Delta}\frac{1 + q^{-1/2}}{\left(1 - q^{-1/2}\right)^2}
$$
where the last equality follows from \eqref{eq-a2}
\end{proof}

\begin{lemma}\label{lem-a6}
Let $\Delta \neq 0$ be some integer.
Then
\begin{equation}\label{eq-lema6}
\sum_{\substack{r_1 \\ (r_1,\Delta) = 1}}
\sum_{\substack{r_2 \\ (r_2,\Delta) = 1}}
\sum_{\substack{d_1 \\ q\mid d_1\Rightarrow q\mid \Delta}}
\sum_{\substack{d_2 \\ q\mid d_2\Rightarrow q\mid \Delta}} \frac{\prod_{\substack{q\mid r_1r_2 \\ q > 3}}\left(1 - \frac{3}{q}\right)^{-1}}{[r_1^2 d_1, r_2^2 d_2]}
\ll 1
\end{equation}
\end{lemma}
\begin{proof}
Denote the sum in \eqref{eq-lema6} by $\Sigma$.
Using the notations from \autoref{lem-a4} we can rewrite $\Sigma$ as 
$$
\Sigma = 
\sum_{\substack{l_r\\(l_r,\Delta) = 1}}
\sum_{\substack{l_d\\q\mid l_d \Rightarrow q\mid \Delta}}
\frac{\tau(l_r)\tau(l_d)}{l_r^2 l_d}\prod_{\substack{q\mid r_1r_2 \\ p > 3}}
\left(1 - \frac{3}{q}\right)^{-1}.
$$
We can write this last sum as an Euler product.
Using \eqref{eq-a2} we get
\begin{multline*}
\Sigma \ll
\prod_{p\nmid 6\Delta}
\left(1 + \left(\frac{1 + p^{-2}}{(1 - p^{-2})^2} - 1\right)\left(1 - \frac{3}{p}\right)^{-1}\right)
\prod_{p\mid \Delta}\left(\frac{1 + p^{-1}}{(1 - p^{-1})^2} \right)\\
\leq c(\Delta) \prod_{p\nmid \Delta}\left(1 + \frac{\BigO{1}}{p^2}\right)
\ll 1.
\end{multline*}
\end{proof}

\subsection{Lemmas for quadratic polynomials}
We now state an important lemma from \cite{Iwaniec-1} regarding representability by a quadratic polynomial $F(x,y)$ in $\mathbb{Z}[x,y]$.
This lemma will help us reduce the problem of reprsentability by $F$ to that of representability by a homogeneous quadratic form.

Let
$$
F(x,y) = ax^2 + bxy + cy^2  + ex + fy + g \in \mathbb{Z}[x,y].
$$
Denote by $\Delta_F$ the discriminant of $F$, and $D_F$ the large discriminant of $F$.
Denote also 
$$
\alpha = bf-2ce,\quad \beta = be-2af.
$$
Assume further that $F$ is in the class $\mathcal{H}$, that $F$ represents arbitrarily large integers, that the partial derivatives of $F$ are linearly independent, and that $\Delta_F$ is not a perfect square.

We will further assume that $F$ satisfies
\begin{equation}\label{eq-ap3}
(a,b,c)\mid (e,f),\quad (ac,\Delta) = (a,b,c)^2, \quad 
(g,\Delta)=1, \quad (\alpha,\beta)\mid \Delta.
\end{equation}

Before stating the lemma, we introduce some notations.
We define $W = (a,b,c,e,f)$ (which is equal to $(a,b,c)$ due to \eqref{eq-ap3}), and $P(x,y) = W^{-1}(F - g)$, which we denote
$$
P(x,y) = a_Px^2 + b_Pxy + c_Py^2 + e_Px + f_Py.
$$
We denote $\Delta_P = W^{-2}\Delta_F$ the discriminant of $P$.
We also denote
$$
\alpha_P = W^{-2}\alpha,\qquad  \beta_P =W^{-2}\beta,
\qquad G_P(x,y) = a_Px^2 + b_Pxy + c_Py^2.
$$

From the equations 
$$
2a_P \alpha_P + b_P \beta_P = \Delta_P e_P,
\qquad
2c_P \beta_P + b_P \alpha_P = \Delta_P f_P,
$$
we have $(2\alpha_P,\Delta_P) = (2\beta_P,\Delta_P) $.
And so, one of $(\alpha_P,\Delta_P), (\beta_P,\Delta_P)$ must divide the other.
Assume without loss of generality that $(\beta_P,\Delta_P) \mid (\alpha_P,\Delta_P)$.
We will denote
$$
\alpha_1 = \frac{\alpha_P }{ (\beta_P,\Delta_P)},
\quad
\beta_1 = \frac{\beta_P }{ (\beta_P,\Delta_P)},
\quad
\Delta_1 = \frac{\Delta_P }{ (\beta_P,\Delta_P)}
$$
and
$$
\varphi(x,y) = G_P\left(\frac{3 + (-1)^{(\alpha_1,\Delta_1)}}{2}x , y\right).
$$

With these notations:
\begin{lemma}[Corollary 2 from \cite{Iwaniec-1}]\label{lem-ap1}
Let $F$ be as above.
Let $A\neq 0$ be some integer.
There exists an integer $Q$ prime to $A$ (which is explicitly given in \cite{Iwaniec-1}) such that the following holds.
If $n$ is an integer such that $n\equiv g (\bmod W)$,
$$
4a_P\frac{n-g}{W} + e_P^2
$$
is a quadratic residue $\bmod(2\beta_P,\Delta_P)$, and
$$
4a_P\frac{n-g}{W} + G_P(\alpha_1,\beta_1) = Q^2\varphi(x,y)
$$
has a solution $x,y$ such that $(x,y,\Delta_1) = 1$
then $n$ is represented by $F$.
\end{lemma}

We now show how this lemma gives a lower bound on $\delta_F(N)$.

\begin{lemma}\label{lem-47}
Let
$$
F(x,y) = ax^2 + bxy + cy^2  + ex + fy + g \in \mathbb{Z}[x,y].
$$
Assume that $F$ represents arbitrarily large integers, that the partial derivatives of $F$ are linearly independent, and that $\Delta_F$ is not a perfect square.
Then
$$
\delta_F(N) \gg \frac{1}{\sqrt{\log N}}
$$
\end{lemma}

\begin{proof}
From \autoref{lem-2} we see that it is enough to prove the statement for the case where $F$ is in $\mathcal{H}$.
We can further assume that $F$ satisfies \eqref{eq-ap3}.
This follows from \cite[Lemma 3]{Iwaniec-1} which states that there is an affine change of variables
$$
F'(x,y) = F(a_1 x + a_2 y + a_3, b_1x + b_2 y + b_3),
$$
such that $F'$ is still in $\mathcal{H}$ and satisfies \eqref{eq-ap3}.
Once more, since $\psi_{F'}(N) \ll \psi_F(N)$, we can simply assume that $F$ satisfies \eqref{eq-ap3}.

Let $Q$ be as in \autoref{lem-ap1} for $A = \Delta D$.
From \cite[Lemma 6]{Iwaniec-1}, there exists an integer $R$ prime to $A$, such that if $n$ is properly represented by the genus of $\varphi$, then $R^2 n$ is represented by $\varphi$ with a representation $R^2 n = \varphi(x_0,y_0)$ which satisfies $(x_0,y_0) \mid R$.
We use the notations from \autoref{lem-ap1}.
Denote 
$$
C = \frac{3 - (-1)^{\Delta_P D_F}}{2}QR.
$$
since $(\alpha_1,\beta_1) = 1$, the number $G_P(\alpha_1,\beta_1)$ is properly represented by $\varphi$.
Let $d$ be the greatest divisor of $G_P(\alpha_1,\beta_1)$ all prime factors of which divide $\Delta_\varphi$.
Let $L$ be defined by the congruence 
$$
C^2 L \equiv \frac{G_P(\alpha_1,\beta_1)}{d} \pmod{\Delta_\varphi},\quad 0<L<\left|\Delta_\varphi\right|.
$$
Denote also 
$$
T = C^2 W \frac{\left| d \Delta_\varphi \right|}{\Delta_1^2},
\quad
l = W \frac{C^2 d L - G_P(\alpha_1,\beta_1)}{\Delta_1^2} + g.
$$
From the proof of \cite[Lemma 13]{Iwaniec-1} we then have that for $n\equiv l \pmod{T}$, the number 
$$
4a_P \frac{n-g}{W} + e_P^2
$$
is a quadratic residue $\bmod (2\beta_P,\Delta_P) $, and $n\equiv g \pmod{W}$.
Thus, from \autoref{lem-ap1} we see that in order for an integer $n$ which is congruent to $l$ mod $T$ to be represented by $F$, it is enough that $\Delta_1^2 (n-g)/W + G(\alpha_1,\beta_1)$ be of the form $Q^2\varphi(x_0,y_0)$ with $x_0,y_0$ satisfying $(x_0,y_0,\Delta_1) = 1$.
We can denote
$$
\frac{\Delta_1^2 (n-g)}{W} + G(\alpha_1,\beta_1) = C^2 d m
$$
with $m\equiv L\, (\bmod{\Delta_\varphi})$.
From our choice of $R$, if $C^2 d m / (QR)^2$ is properly represented by the genus of $\varphi$, then $C^2 d m$ is of the form $Q^2 \varphi(x_0,y_0)$ with $x_0,y_0$ satisfying $(x_0,y_0)\mid R$.
Since  $(R,\Delta_1) = 1$, we also have $(x_0,y_0,\Delta_1) = 1$.
From \cite[Lemma 7]{Iwaniec-1}, the condition $q\mid m\Rightarrow \left(\frac{\Delta_\varphi}{q}\right) = 1$ is sufficient for the number $C^2 d m / (QR)^2$ to be properly represented by the genus of $\varphi$.

Thus,
$$
\sum_{\substack{0<n \leq N \\ n = F(x,y)}}1 
\gg
\sum_{\substack{m\leq N \\ m\equiv L (\bmod{\Delta_\varphi}) \\ q\mid m \Rightarrow \left(\frac{\Delta_\varphi}{q}\right) = 1}} 1.
$$

The last sum can be bounded from below using a $\frac{1}{2}$-dimensional sieve.
We start with the set 
$\mathcal{M} = \left\lbrace m \;\middle| \; m\leq N, \; m\equiv L \pmod{\Delta_\varphi} \right\rbrace$.
For $\lambda$ prime to $\Delta_\varphi$, all prime factors of which satisfy $\left(\frac{\Delta_\varphi}{q}\right) = -1$ we denote $\mathcal{M}_\lambda = \left\lbrace m\in \mathcal{M} \; \middle| \; \lambda \mid m \right\rbrace$.
Then $\left|\mathcal{M}_\lambda\right| = \frac{N}{\lambda \left|\Delta_\varphi\right|} + \mathcal{O}\left(1\right)$.

Since each $m\in\mathcal{M}$ satisfies $\left(\frac{\Delta_\varphi}{m}\right) = 1$, we see that each $m$ has to be divisible by an even number of primes $q$ for which  $\left(\frac{\Delta_\varphi}{q}\right) = -1$.
It follows that it is enough to sieve out primes up to $\sqrt{N}$.

And so, applying \autoref{thm-sieves} with $D = N (\log N)^{-2}$ and some $1<s<2$ we get
\begin{multline*}
\sum_{\substack{m\leq N \\ m\equiv L (\bmod{\Delta_\varphi}) \\ q\mid m \Rightarrow \left(\frac{\Delta_\varphi}{q}\right) = 1}} 1
\gg \\
N \prod_{\substack{p \leq D^{1/s} \\ \left(\frac{\Delta_\varphi}{p}\right) = -1}}\left(1 - \frac{1}{p}\right)\left\lbrace c_{\frac{1}{2}}(s) + \mathcal{O}\left(\log^{-1/6} (D) \right) \right\rbrace + \mathcal{O}\left(\frac{N}{\log^2 N}\right)
\\ \gg 
\frac{N}{\sqrt{\log N}}.
\end{multline*}
And so
$$
\delta_F(N) \gg 
\frac{1}{N}
\sum_{\substack{m\leq N \\ m\equiv L (\bmod{\Delta_\varphi}) \\ q\mid m \Rightarrow \left(\frac{\Delta_\varphi}{q}\right) = 1}} 1
\gg 
\frac{1}{\sqrt{\log N}}.
$$
\end{proof}
\begin{lemma}\label{lem-48}
Let $F$ be as in \autoref{lem-47}.
Then
$$
\delta_F(N) \ll \frac{1}{\sqrt{\log N}}.
$$
\end{lemma}
\begin{proof}
Denote $G(x,y) = ax^2 + bxy + cy^2$.
One can check that
$$
\Delta_F^2 F(x,y) = G(\Delta_F x + \alpha, \Delta_F y + \beta) + D_F \Delta_F.
$$
Denote $G'(x,y) = (a,b,c)^{-1}G(x,y)$.
It follows that $\delta_F(N) \ll \delta_{G'}(N)$.
From genus theory for quadratic forms, it is known that if $n = G'(x,y)$ then $n$ must be of the form $n = r^2 d m$ where $q\mid d \Rightarrow q\mid \Delta_{G'}$, and $q\mid m \Rightarrow \left(\frac{\Delta_{G'}}{q}\right) = 1$.

Thus, we have that 
\begin{equation}\label{eq-49}
\sum_{\substack{n\leq N \\ n = F(x,y)}}1
\ll
\sum_{r}\sum_{\substack{d \\ q\mid d \Rightarrow q\mid \Delta_{G'}}}
\sum_{\substack{m\leq N/r^2 d \\ q\mid m \Rightarrow \left(\frac{\Delta_{G'}}{q}\right) = 1}} 1.
\end{equation}
Once more, we can bound the last sum using the $\frac{1}{2}$-dimensional sieve.
for $X>0$, let $\mathcal{M} = \left\lbrace m \; \middle| \; m\leq X \right\rbrace $, and $\mathcal{M}_\lambda = \left\lbrace m\in \mathcal{M} \; \middle| \; \lambda \mid m \right\rbrace$.
Let 
$$
\mathcal{P} = \left\lbrace p \; \middle| \; \left(\frac{\Delta_{G'}}{p}\right) = -1 \right\rbrace.
$$
Applying theorem \autoref{thm-sieves} with $D = X/\log^2 X$ and some $s\geq 1$ we get that
\begin{multline*}
\sum_{\substack{m\leq X \\ q\mid m \Rightarrow \left(\frac{\Delta_{G'}}{q}\right) = 1}} 1
\ll\\
X \prod_{\substack{p \leq D^{1/s} \\ p\in \mathcal{P}}}\left(1 - \frac{1}{p}\right)\left\lbrace C_{\frac{1}{2}}(s) + \mathcal{O}\left(\log^{-1/6} (D) \right) \right\rbrace + \mathcal{O}\left(\frac{X}{\log^2 X}\right)\\
\ll
\frac{X}{\sqrt{\log (X + 3)}}
\end{multline*}

Thus, plugging this into \eqref{eq-49} we get that
$$
\sum_{\substack{n\leq N \\ n = F(x,y)}}1
\ll
\sum_{r}\sum_{\substack{d \\ q\mid d \Rightarrow q\mid \Delta_{G'}}}
\frac{N}{r^2 d \sqrt{\log\left(\frac{N}{r^2 d} + 3 \right)}}
\ll \frac{N}{\sqrt{\log N}}.
$$
It follows that $\delta_F(N) \ll \frac{1}{\sqrt{\log N}}$.
\end{proof}

\section{Small multiples of large primes represented by F}\label{sec-4}
In this section we prove \autoref{lem-43}, \autoref{lem-44} and \autoref{lem-45}.
These lemmas all have to do with estimating the number of small multiples of large primes which are represented by $F$.
Throughout the section, we let $F$ be as in \autoref{main-thm3}, and we use the same notations.

We further denote
$$
\alpha = bf -2ce,\quad \beta=be - 2af,
$$
and $G_F(x,y) = ax^2 + bxy + cy^2$ the quadratic part of $F$.

\subsection{an upper bound for $S_k(N)$}
We now prove \autoref{lem-43}.
\begin{proof}[proof of \autoref{lem-43}]
From computation, one can see that
$$
\Delta^2 F(x,y) = G_F(\Delta x + \alpha, \Delta y + \beta ) + \Delta D.
$$
And so,
\begin{multline}\label{eq-5-1}
S_k(N) =\\
\sum_{\substack{p \text{ prime} \\ kp = F(x,y) \leq N}} 1
\ll 
\sum_{\substack{p \text{ prime} \\ kp \leq N \\ \Delta^2 k p = G_F(\Delta x + \alpha, \Delta y + \beta ) + \Delta D }} 1 \ll 
\sum_{\substack{p \text{ prime} \\ kp \leq N \\ \Delta^2 k p - \Delta D = G_F( x ,  y  )   }} 1 .
\end{multline}

We denote $C = (a,b,c)$, $B = -\Delta D$, $A = \Delta^2 k$.
Denote also $\varphi = C^{-1}G_F$, the primitive quadratic form associated to $G_F$.
Note that $\Delta_\varphi = C^{-2} \Delta$ which means that 
$\left(\frac{\Delta_\varphi}{q}\right) = \left(\frac{\Delta}{q}\right)$ for all primes not dividing $C$.

Throughout the proof, the implied constants in the $\ll$ notation will depend on $F$, but not on $k$. 
This means that they can depend on $B,C,\Delta$, but not on $A$.

With these notations we see from \eqref{eq-5-1} that
\begin{equation}\label{eq-50}
S_k(N) \ll \sum_{\substack{p \text{ prime} \\ p \leq N/A \\ Ap + B = C\varphi(x,y)}}1.
\end{equation}

This sum was considered by Iwaniec in \cite[Theorem 3]{Iwaniec-1}.
We follow the proof of the upper bound given there, with a slight modification so that the bound remains uniform for $A \leq \log^5 N$.

From the proof of \cite[Theorem 3]{Iwaniec-1} we get
\begin{multline}\label{eq-51}
\sum_{\substack{p \leq N/A \\ Ap + B = C\varphi(x,y)}}1
\leq  \\
\sum_{\substack{d \\ q\mid d \Rightarrow q\mid\Delta}}
\sum_{\substack{r \\ 
\left( \frac{AB}{(A,B)^2}, \frac{Cr^2 d}{(A,B,Cr^2 d)}\right) = 1}}
\sum_{\substack{m  \leq N  + B \\ 
\left( m, \frac{AB}{(A,B)^2}\right) = 1 \\ 
q\mid m \Rightarrow \left(\frac{\Delta}{q}\right) = 1}}
\sum_{\substack{x \leq \frac{N(A,B,Cr^2 d)}{A Cr^2 d m} \\
|a_ix + b_i| \text{ primes for }i=1,2}} 1
\end{multline}
where 
\begin{equation*}
a_1 = \frac{C r^2 d m}{(A,B,Cr^2 d )}, \; \; 
b_1 = l_{m}, \; \;
a_2 = \frac{A}{(A,B)}, \;\;
b_2 = \frac{(Al_{m} + B)(A,B,Cr^2 d)}{(A,B)Cr^2 d m}
\end{equation*}
and $l_{m}$ is a solution to the congruence
$$
\frac{A}{(A,B)}l_{m} + \frac{B}{(A,B)} \equiv 0\left( \bmod \frac{Cr^2 d m}{(A,B,Cr^2 d )}\right).
$$

Denote
$$
\Sigma(r,d,m) = \sum_{\substack{x \leq \frac{N(A,B,Cr^2 d)}{A Cr^2 d m} \\
|a_ix + b_i| \text{ primes for }i=1,2}} 1.
$$
We get from \cite[Lemma 9]{Iwaniec-1} that
\begin{multline}\label{eq-51.5}
\Sigma(r,d,m) \ll \\
\frac{N}{Adr^2 m} 
\log^{-2}\left(\frac{N}{Adr^2 m}  + 3 \right)
\prod_{\substack{q\mid r m}}\left(1 - \frac{1}{q}\right)^{-2}
\prod_{\substack{q\mid A}}\left(1 - \frac{1}{q}\right)^{-1}.
\end{multline}

By \eqref{eq-51.5}, \autoref{lem-a3}, and using partial summation we get that

\begin{multline*}
\sum_{\substack{m  \leq N  + B \\ 
\left( m, \frac{AB}{(A,B)^2}\right) = 1 \\ 
q\mid m \Rightarrow \left(\frac{\Delta}{q}\right) = 1}} \Sigma(r,d,m) 
\ll 
\frac{N}{Adr^2}
\prod_{\substack{q\mid r}}\left(1 - \frac{1}{q}\right)^{-2}
\prod_{\substack{q\mid A}}\left(1 - \frac{1}{q}\right)^{-1}
\times \\ \times
\sum_{\substack{m  \leq N  + B \\ 
\left( m, \frac{AB}{(A,B)^2}\right) = 1 \\ 
q\mid m \Rightarrow \left(\frac{\Delta}{q}\right) = 1}}
\frac{1}{m}\log^{-2}\left(\frac{N}{Adr^2 m}  + 3 \right)
 \ll \\
\frac{N}{A dr^2 \log^{3/2} N} 
\prod_{\substack{q\mid r}}\left(1 - \frac{1}{q}\right)^{-2}
\prod_{\substack{q \mid A \\ \left(\frac{\Delta }{q}\right) = -1}}\left( 1 - \frac{1}{q}\right)^{-1}.
\end{multline*}
Plugging this into \eqref{eq-51} we get
\begin{multline}\label{eq-52}
\sum_{\substack{p \leq N/A \\ Ap + B = C\varphi(x,y)}}1
\ll \\  
\sum_{\substack{d \\ q\mid d \Rightarrow q\mid\Delta}}
\sum_r
\frac{N }{A dr^2 \log^{3/2} N} 
\prod_{\substack{q\mid r}}\left(1 - \frac{1}{q}\right)^{-2}
\prod_{\substack{q \mid A \\ \left(\frac{\Delta }{q}\right) = -1}}\left( 1 - \frac{1}{q}\right)^{-1}
=\\
\frac{N}{A\log^{3/2} N}\prod_{\substack{q \mid A \\ \left(\frac{\Delta }{q}\right) = -1}}\left( 1 - \frac{1}{q}\right)^{-1}
\sum_{r}\left(\frac{1}{r^2}
\prod_{\substack{q\mid r}}\left(1 - \frac{1}{q}\right)^{-2}\right)
\sum_{\substack{d \\ q\mid d \Rightarrow q\mid\Delta}} \frac{1}{d}
.
\end{multline}

Noting that 
$$
\sum_r \frac{1}{r^2}\prod_{q\mid r}\left(1 - \frac{1}{q}\right)^{-2},\qquad 
\sum_{\substack{d \\ q\mid d \Rightarrow q\mid \Delta}}\frac{1}{d}
= \prod_{q\mid \Delta}\left(1 - \frac{1}{q}\right)^{-1}
$$
are finite, and also noting that $A = \Delta^2 k$, we get from \eqref{eq-52} that
$$
\sum_{\substack{p \leq N/A \\ Ap + B = C\varphi(x,y)}}1
\ll
\frac{N}{k\log^{3/2}N} 
\prod_{\substack{q \mid k \\ \left(\frac{\Delta }{q}\right) = -1}}\left( 1 - \frac{1}{q}\right)^{-1}.
$$

From \eqref{eq-50} we then get the required bound
$$
S_k(N) \ll \frac{N}{k\log^{3/2}N} 
\prod_{\substack{q \mid k \\ \left(\frac{\Delta }{q}\right) = -1}}\left( 1 - \frac{1}{q}\right)^{-1}.
$$
\end{proof}
\subsection{A lower bound for $S_k(N)$}
We now prove \autoref{lem-44}.
In \cite[Lemma 13]{Iwaniec-1} Iwaniec proved a lower bound for $S_1(N)$.
We follow the proof given there, with slight modifications so that we can prove lower bounds for all $S_k(N)$ with $k\leq \log^5 N$.

We can assume that our function $F$ is in the class $\mathcal{H}$. 
This is true due to \autoref{lem-2}, replacing $F$ by the function $H$ given in the lemma if necessary.

We can further assume that $F$ satisfies
\begin{equation}\label{eq-53}
(a,b,c)\mid (e,f),\quad (ac,\Delta) = (a,b,c)^2, \quad 
(g,\Delta)=1, \quad (\alpha,\beta)\mid \Delta.
\end{equation}
This follows from \cite[Lemma 3]{Iwaniec-1} which states that there is an affine change of variables
$$
F'(x,y) = F(a_1 x + a_2 y + a_3, b_1x + b_2 y + b_3),
$$
such that $F'$ is still in $\mathcal{H}$ and satisfies \eqref{eq-53}.
Once more, since $\psi_{F'}(N) \ll \psi_F(N)$, we can simply assume that $F$ satisfies \eqref{eq-53}.

We now prove \autoref{lem-44}.
\begin{proof}[Proof of \autoref{lem-44}]
We use the notations of \autoref{lem-ap1}.
Let $Q$ be as \autoref{lem-ap1} for $A = \Delta D$.
From \cite[Lemma 6]{Iwaniec-1}, there exists an integer $R$ prime to $A$, such that if $n$ is properly represented by the genus of $\varphi$, then $R^2 n$ is represented by $\varphi$ with a representation $R^2 n = \varphi(x_0,y_0)$ which satisfies $(x_0,y_0) \mid R$.

Denote 
$$
C = \frac{3 + (-1)^{\Delta_P D_F}}{2}QR, \quad
$$
and $\Delta_\varphi$ the discriminant of $\varphi$.
Let $d$ be the largest divisor of $G_P(\alpha_1, \beta_1)$ all prime factors of which divide $\Delta_\varphi$.
Let $L$ be defined by the congruence
$$
C^2 L \equiv \frac{G_P(\alpha_1, \beta_1)}{d} \pmod{\Delta_\varphi} ,\quad 0<L<\left|\Delta_\varphi\right|.
$$

Finally, let 
$$
T = C^2 W \frac{\left| d \Delta_\varphi \right|}{\Delta_1^2},
\quad
l = W \frac{C^2 d L - G_P(\alpha_1,\beta_1)}{\Delta_1^2} + g.
$$
From the proof of \cite[Lemma 13]{Iwaniec-1} we then have that for $n\equiv l \pmod{T}$, the number 
$$
4a_P \frac{n-g}{W} + e_P^2
$$
is a quadratic residue $\bmod (2\beta_P,\Delta_P) $, and $n\equiv g \pmod{W}$.

We choose $c_F = T\Delta D$, that is we assume that $(k,T) = 1$ and also $(k, \Delta D) = 1$.
Denote
$$
\mathcal{M} = \left\lbrace
m \,\middle|\, \Delta_1^2 \frac{kp - g}{W} + G_P(\alpha_1,\beta_1) = C^2 dm,\; kp\equiv l \pmod{T}, \; kp \leq N
\right\rbrace.
$$
It follows from \autoref{lem-ap1} and from our choice of $R$ that for $m\in\mathcal{M}$, if $C^2 d m / R^2$ is properly represented by the genus of $\varphi$, then $kp$ is represented by $F$.
And so, we have reduced the problem of representability by $F$ to that of proper representability by the genus of $\varphi$.

Genus theory for quadratic forms gives simple criteria for when a number is properly represented by the genus of $\varphi$.
Such a lemma is given in \cite[Lemma 5]{Iwaniec-1}.
From which we get that for $C^2 d m / R^2$ to be properly represented by the genus of $\varphi$, it suffices that all prime divisors $q$ of $m$ satisfy $\left(\frac{\Delta_\varphi}{q}\right) = 1$.

It can be verified that for $m\in\mathcal{M}$,
$$
(m, \Delta_P D W^{-2}) = 1, \quad (m,k) = 1.
$$
Thus, denoting 
$$
\mathcal{Q} = \left\lbrace
q \,\middle|\, \left(\frac{\Delta_\varphi}{q}\right)=1;\quad (q,k\Delta_P D / W^2)=1\right\rbrace,
$$
we get that
\begin{equation}\label{eq-53.25}
S_k(N) \gg \sum_{\substack{m \in \mathcal{M} \\ q\mid m \Rightarrow q\in\mathcal{Q}}} 1.
\end{equation}

In order to estimate this sum, we will use the $\frac{1}{2}$-dimensional sieve.
For $\lambda$ such that $(\lambda , k\Delta_P D / W^2) = 1$ denote
\begin{multline*}
\mathcal{M}_\lambda
= \left\lbrace m \in \mathcal{M} \,\middle|\, \lambda \mid m\right\rbrace
= \\
\left\lbrace m \,\middle|\, \Delta_1^2 \frac{kp - g}{W} + G_P(\alpha_1,\beta_1) = C^2 d m, \; kp\equiv l_\lambda \pmod{\lambda T},\; kp \leq N\right\rbrace
\end{multline*}
where
\begin{equation*}
l_\lambda = W \frac{C^2 d L_\lambda - G_P(\alpha_1,\beta_1)}{\Delta_1^2} + g,\qquad L_\lambda \equiv 
\begin{cases}
L \pmod{\Delta_\varphi}\\
0 \pmod{\lambda}
\end{cases}.
\end{equation*}
We can then see that $\left|\mathcal{M}_\lambda\right| = \pi\left(N/k, \lambda T, l_\lambda\right) + \mathcal{O}(1)$.
Denote
$$
\mathcal{P} = \left\lbrace
q \,\middle|\, \left(\frac{\Delta_\varphi}{q}\right)=-1;\quad (q,k\Delta_P D / W^2)=1\right\rbrace.
$$
All $m\in\mathcal{M}$ satisfy $m\equiv L (\bmod \Delta_\varphi)$.
From \cite[Lemma 4]{Iwaniec-1} we have that $\left(\frac{\Delta_\varphi}{L}\right) = 1$.
It follows that $m\in\mathcal{M}$ satisfy $\left(\frac{\Delta_\varphi}{m}\right) = 1$.
And so, each $m\in\mathcal{M}$ must be divisible by an even number of primes from $\mathcal{P}$.

Denote by $\mathcal{S}(\mathcal{A},z)$ the size of the sifted set 
$$
\mathcal{A}\setminus\bigcup_{\substack{p\in\mathcal{P} \\ p< z}}\mathcal{A}_p.
$$
Since each $m\in\mathcal{M}$ is divisible by an even number of primes from $\mathcal{P}$, we have that for $1\leq s < 2$
\begin{equation}\label{eq-53.5}
\sum_{\substack{m\in\mathcal{M} \\ q\mid m \Rightarrow q\in\mathcal{Q}}} 1
\gg \; \mathcal{S}(\mathcal{M},N^{1/2s})
-
\sum_{\substack{N^{1/2s}\leq p_1 \leq p_2 \\ p_1,p_2\in\mathcal{P}}}
\mathcal{S}(\mathcal{M}_{p_1p_2},N^{1/2s}).
\end{equation}

The term $\mathcal{S}(\mathcal{M},N^{1/2s})$ can be evaluated using the $\frac{1}{2}$-dimensional sieve.
We apply \autoref{thm-sieves}, where we set
$\kappa = 1/2$, $D = \sqrt{N} / \log^{40} N$, $z = D^{1/s}\sim N^{1/2s}$, $X = \text{Li}(N/k)/\phi(T)$ and $g(p) = \phi(T)/\phi(pT)$ for $p\in\mathcal{P}$ and $g(p)=0$ for primes not in $\mathcal{P}$.
We then get
\begin{multline}\label{eq-53.6}
\mathcal{S}(\mathcal{M},N^{1/2s}) \geq\\
\frac{N}{k\log N}
\prod_{\substack{p<z\\p\in\mathcal{P}}}
\left(1 - \frac{\phi(T)}{\phi(pT)}\right)
\left\lbrace c_{1/2}(s) + \mathcal{O}\left(\log^{-1/12} N\right) \right\rbrace + \mathcal{R}(D)
\end{multline}
where
$$
\mathcal{R}(D) = 
\sum_{\substack{\lambda \leq D \\ q\mid \lambda \Rightarrow q\in\mathcal{P}}}\left|\pi\left(\frac{N}{k},\lambda T, \frac{l_\lambda}{k}\right) - \frac{\text{Li}(N/k)}{\phi(\lambda T)}\right|.
$$

Using the \nameref{b-v}, and since $k\leq \log^5 N$, we get that $\mathcal{R}(D) \ll N\log^{-30} N$ which is going to be negligible.
And so, since we regard $F$ as a constant, we get from \eqref{eq-53.6} that
\begin{multline}\label{eq-54}
\mathcal{S}(\mathcal{M},N^{1/2s}) \gg \\
\frac{N}{k\log^{3/2} N}\prod_{\substack{q\mid k \\ \left(\frac{\Delta}{q}\right) = -1}}\left(1 - \frac{1}{q}\right)^{-1}\left\lbrace c_{1/2}(s)
+\mathcal{O}\left(\frac{1}{\log^{1/12} N}\right)\right\rbrace.
\end{multline}
We further note that $c_{1/2}(s)$ behaves like $\sqrt{e^\gamma / \pi} (s-1)^{1/2}$ near $s = 1$.

As for the term 
$$
\sum_{\substack{N^{1/2s}\leq p_1 \leq p_2 \\ p_1,p_2\in\mathcal{P}}}
\mathcal{S}(\mathcal{M}_{p_1p_2},N^{1/2s})
$$
we have that
\begin{multline}\label{eq-55}
\sum_{\substack{N^{1/2s}\leq p_1 \leq p_2 \\ p_1,p_2\in\mathcal{P}}}
\mathcal{S}(\mathcal{M}_{p_1p_2},N^{1/2s}) 
\leq  \\
\sum_{\substack{m\leq M N ^{-1/s} \\ q\mid m \Rightarrow q\in\mathcal{Q}}} \quad
\sum_{\substack{N^{1/2s}\leq p_1 \leq \sqrt{M} \\ p_1 \in \mathcal{P}}}
\sum_{\substack{kp\leq N \\ kp\equiv l_{p_1 m} \pmod{p_1 m T} \\ \Delta_1^2(kp - g) + W G_P(\alpha_1,\beta_1)=WC^2dp_1p_2m \\ p_1, p_2 \text{ prime}}}1 \\
=\sum_m \sum_{p_1} \Sigma(m,p_1)
\end{multline}
where 
$$
M = \frac{\Delta_1^2(N -g) + WG_P(\alpha_1,\beta_1)}{WC^2d}.
$$

The inner sum $\Sigma$ can be written as
$$
\Sigma(m,p_1) = 
\sum_{\substack{x\leq N/kp_1 mT \\ 
p_1 m T x + (k^{-1} l_{p_1 m} (\bmod p_1 m T)) \text{ prime}\\
k|\Delta_\varphi|x + \frac{\Delta_1^2l_{p_1m} + WG_P(\alpha_1,\beta_1) - \Delta_1^2g}{WC^2dp_1m} \text{ prime}
}} 1.
$$

And so, from \cite[Lemma 9]{Iwaniec-1} we get
$$
\Sigma(m,p_1) \ll
\frac{N}{kp_1 m T \log^2(N / k p_1 m T)}
\prod_{\substack{q\mid m }}\left(1 - \frac{1}{q}\right)^{-1}
\prod_{\substack{q\mid k}}\left(1 - \frac{1}{q}\right)^{-1}.
$$

Thus, from \eqref{eq-55} we get
\begin{multline}\label{eq-56}
\sum_{\substack{N^{1/2s}\leq p_1 \leq p_2 \\ p_1,p_2\in\mathcal{P}}}
\mathcal{S}(\mathcal{M}_{p_1p_2},N^{1/2s}) 
\ll  \\
\frac{N}{k\log^2 N}\prod_{q\mid k}\left(1-\frac{1}{q}\right)^{-1}
\sum_{\substack{m \leq N^{1 - 1/s} \\ q\mid m \Rightarrow q\in\mathcal{Q}}} \frac{1}{m}\prod_{q\mid m}\left(1-\frac{1}{q}\right)^{-1}
\sum_{\substack{N^{1/2s} < p < N^{1/2} \\ p \in \mathcal{P}}}\frac{1}{p}\\
= 
\frac{N}{k\log^2 N}\prod_{q\mid k}\left(1-\frac{1}{q}\right)^{-1}
\times \Sigma_m \times  \Sigma_p.
\end{multline}

For $\Sigma_m$ we have
$$
\Sigma_m \ll \prod_{\substack{q\leq N^{1- 1/s} \\ (q,k)=1 \\ \left(\frac{\Delta}{q}\right) = 1}}\left(1  + \frac{1}{q-1}\right)
\ll \sqrt{1-\frac{1}{s}}\log^{1/2} N 
\prod_{\substack{q\mid k \\ \left(\frac{\Delta}{q}\right) = 1}}\left(1-\frac{1}{q}\right).
$$
As for $\Sigma_p$ we have $\Sigma_p \ll \log s$.

Plugging these results into \eqref{eq-56} we get
\begin{equation}\label{eq-57}
\sum_{\substack{N^{1/2s}\leq p_1 \leq p_2 \\ p_1,p_2\in\mathcal{P}}}
\mathcal{S}(\mathcal{M}_{p_1p_2},N^{1/2s}) 
\ll
C(s)
\frac{N}{k\log^{3/2} N}
\prod_{\substack{q\mid k \\ \left(\frac{\Delta}{q}\right) = -1}}\left(1-\frac{1}{q}\right)^{-1}
\end{equation}

where $C(s) = \log s\sqrt{(s-1)/s}$ behaves like $(s-1)^{3/2}$ near $s=1$.

And so, from \eqref{eq-53.25},\eqref{eq-53.5},\eqref{eq-54}, \eqref{eq-57} we get that
$$
S_k(N) \geq \left(K_1 c_{1/2}(s) -K_2 C(s)\right)\frac{N}{k\log^{3/2} N}
\prod_{\substack{q\mid k \\ \left(\frac{\Delta}{q}\right) = -1}}\left(1-\frac{1}{q}\right)^{-1}
$$
for some positive constants $K_1,K_2$ which depend on $F$.
From the behaviour of $c_{1/2}(s),C(s)$ near $s=1$, we get that for $s$ sufficiently close to $1$, $K_1 c_{1/2}(s) > K_2 C(s)$.
And so
$$
S_k(N) \gg \frac{N}{k\log^{3/2} N}
\prod_{\substack{q\mid k \\ \left(\frac{\Delta}{q}\right) = -1}}\left(1-\frac{1}{q}\right)^{-1}
$$
as required.
\end{proof}

\subsection{An upper bound for $S_{k_1,k_2}(N)$}
We now prove \autoref{lem-45}.
\begin{proof}[Proof of \autoref{lem-45}]
Let $k_1<k_2\leq \log^5 N$.
Using the identity
$$
\Delta^2 F(x,y) = G(\Delta x + \alpha, \Delta y + \beta) + \Delta D
$$
we get that
$$
S_{k_1,k_2}(N)\ll
\sum_{\substack{p\leq N/k_1 \text{ prime}\\ \Delta^2 k_1 p = G(x_1,y_1) + \Delta D \\ \Delta^2 k_2p =  G(x_2,y_2) + \Delta D}} 1.
$$
From genus theory for quadratic forms (e.g. \cite[Lemma 5]{Iwaniec-1}), we get that if $n$ is properly represented by the genus of $G$, then every prime divisor $q$ of $n$ which is prime to $\Delta$ must satisfy $\left(\frac{\Delta}{q}\right) = 1$.
Thus, we get that
\begin{multline}\label{eq-58}
S_{k_1,k_2}(N) \ll \\
\sum_{\substack{r_1,r_2\in\mathbb{N}\\ (r_1r_2,\Delta) = 1}} 
\sum_{\substack{d_1,d_2 \\ q\mid d_1 d_2 \Rightarrow q \mid \Delta }}
\sum_{\substack{m_1,m_2 \\ (m_1 m_2,\Delta) = 1 \\ q\mid m_1m_2 \Rightarrow \left(\frac{\Delta}{q}\right) = 1 \\ m_1,m_2 \text{ square free}}}\;
\sum_{\substack{p\leq N/k_1 \text{ prime}\\ \Delta^2 k_1 p - \Delta D = r_1^2 d_1 m_1  \\ \Delta^2 k_2 p - \Delta D = r_2^2 d_2 m_2}} 1 \\
= \sum_{r_1,r_2}\;\sum_{d_1,d_2}\;\Sigma_{m_1,m_2,p}.
\end{multline}

We now wish to give an upper bound for the last sum.
For this, we split the sum into several cases.
We fix some small $\epsilon>0$.
\begin{case}
We first consider the case where $[r_1^2 d_1, r_2^2 d_2]\leq N^\epsilon$.
In this case we give an upper bound for $\Sigma_{m_1,m_2,p}$ using a $2$ dimensional sieve.
Note that for some choices of $r_1,r_2,d_1,d_2$ we will have $\Sigma_{m_1,m_2,p}=0$.
In this case $\Sigma_{m_1,m_2,p}$ will satisfy the upper bound trivially.
And so, in the argument that follows we assume that there are no local obstructions to the equations
\begin{equation}\label{eq-511}
\Delta^2k_i p - \Delta D = r_i^2d_im_i,\;\text{i=1,2}.
\end{equation}

In order for the equation $\Delta^2 k_1 p - \Delta D = r_1^2d_1m_1$ to have a solution, we must have $\Delta(\Delta k_1, D)\mid r_1^2 d_1 m$.
We then have
$$
\frac{\Delta k_1}{(\Delta k_1, D)}p
-\frac{D}{(\Delta k_1, D)} = 
\frac{r_1^2 d_1 }{\Delta (\Delta k_1,D,r_1^2 d_1)} m_1'
$$
where 
$$
m_1' = \frac{(\Delta k_1,D,r_1^2 d_1)}{(\Delta k_1, D)} m_1.
$$
Denote $T_1 = \frac{r_1^2 d_1 }{\Delta (\Delta k_1,D,r_1^2 d)} $ and
$$
l_1  = \frac{D}{(\Delta k_1, D)}\left(\frac{\Delta k_1}{(\Delta k_1, D)}\right)^{-1}\pmod{T_1}.
$$
It follows that we must have
$$
p\equiv l_1 \pmod{T_1}.
$$
We can similarly define $T_2,l_2$ and get that in order for
$$
\Delta^2 k_2 p - \Delta D = r_2^2d_2m_2
$$
to have a solution, we must have $p \equiv l_2 \pmod{T_2}$.
There exists $l$ which satisfies
$$
l \equiv l_1 \pmod{T_1},\quad l \equiv l_2 \pmod{T_1}
$$
for otherwise we would have a local obstruction.

Let 
$$
\mathcal{A} = \left\lbrace n\leq N/k_1 \; \middle| \; n\equiv l \pmod{[T_1,T_2]}\right\rbrace.
$$

We would now like to sift $\mathcal{A}$ for primes $p$, which also satisfy \eqref{eq-511}.
To this end, we define for each prime $q$ a set $\Omega_q$ of residue classes mod $q$ that $n\in\mathcal{A}$ should avoid.
For simplicity, we choose not to sift by primes dividing $6\Delta D k_1 k_2 (k_1 - k_2) r_1 r_2$.
Denote the set of these primes by $\mathcal{Q}$.

For $q$ such that $\left(\frac{\Delta}{q}\right) = -1$ and $q\not \in \mathcal{Q}$ we define 
\begin{multline*}
\Omega_q = \left\lbrace 0, 
\frac{D}{(\Delta k_1, D)}\left(\frac{\Delta k_1}{(\Delta k_1, D)}\right)^{-1}\pmod{q},\right. \\ \left.
\frac{D}{(\Delta k_2, D)}\left(\frac{\Delta k_2}{(\Delta k_2, D)}\right)^{-1}\pmod{q}\right\rbrace.
\end{multline*}
For other $q\not \in \mathcal{Q}$, we define $\Omega_q =  \left\lbrace 0\right\rbrace$.
For $q\in\mathcal{Q}$ we define $\Omega_q = \phi$.
Define also 
\begin{equation*}
g(q) = \frac{\left|\Omega_q\right|}{q} = 
\begin{cases}
1/q & \left(\frac{\Delta}{q}\right) = 1 \text{ and } q\not \in \mathcal{Q} \\
3/q & \left(\frac{\Delta}{q}\right) = -1 \text{ and } q\not \in  \mathcal{Q} \\ 
0 & q\in\mathcal{Q}.
\end{cases}.
\end{equation*}
We will now use the fact that for $z < N$,
$$
\Sigma_{m_1,m_2,p} \ll
\mathcal{S}(\mathcal{A},z,\Omega) + \BigO{z},
$$
and we will apply \autoref{thm:sieves-variant} in order to get an upper bound for $\mathcal{S}(\mathcal{A},z,\Omega)$.

In \autoref{thm:sieves-variant} we set $X = \frac{N/k_1}{[T_1,T_2]}$ and $g$ as above (where we extend $g$ to be multiplicative).
With these choices, we see that $r_\lambda = |\mathcal{A}^\Omega_\lambda | - Xg(\lambda) $ satisfies $|r_\lambda| \leq g(\lambda)\lambda$ for square-free $\lambda$'s.
From the definition of $g$, we see that this implies the bound $|r_\lambda| \leq d_3(\lambda)$ on square-free $\lambda$'s (where $d_3$ is the 3-fold divisor function).
It follows that 
$$
\sum_{\substack{\lambda < D \\ \lambda \text{ square free}}}\left|r_\lambda \right| \ll D\log^2 D.
$$
And so, in \autoref{thm:sieves-variant} we will set $D = N^{1-\epsilon}\log^{-100}N$, which will make the error term $\sum_{\substack{\lambda < D \\ \lambda \mid P(z) }}|r_\lambda| $ negligible.
Note that from the definition of $T_1,T_2$ we have that $T_i \mid r_i^2 d_i$. 
It follows that $[T_1,T_2]\ll [r_1^2d_1,r_2^2d_2]$ which implied
$$
\frac{N}{[T_1,T_2]} \gg \frac{N}{[r_1^2d_1,r_2^2d_2]} \gg N^{1 - \epsilon}.
$$

We now apply \autoref{thm:sieves-variant} with $\kappa = 2$ and $z = D^{1/s}$.
We choose $s$ suitably large such that $C_2(s) > 0$.
We then get that 
\begin{multline*}
\Sigma_{m_1,m_2,p} \ll \mathcal{S}(\mathcal{A}, z, \Omega) + \BigO{z} \ll \\
\frac{N/k_1}{[r_1^2 d_1 , r_2^2 d_2]\log^2 N} \prod_{\substack{q\in\mathcal{Q} \\ q > 3}}\left( 1- \frac{3}{q}\right)^{-1}
+ \BigO{ N^{1 - \epsilon}\log^{-98}N}\ll \\
\frac{N/k_1}{[r_1^2 d_1 , r_2^2 d_2]\log^2 N} 
\prod_{ \substack{q\mid k_1 k_2 (k_1 - k_2)\\ q > 3} }\left( 1- \frac{3}{q}\right)^{-1}
\prod_{ \substack{q\mid r_1r_2\\ q > 3} }\left( 1- \frac{3}{q}\right)^{-1}.
\end{multline*}
Denote 
$$
C(k_1,k_2) = \prod_{ \substack{q\mid k_1 k_2 (k_1 - k_2) \\ q > 3 }}\left( 1- \frac{3}{q}\right)^{-1}.
$$
So in the case $[r_1^2d_1,r_2^2 d_2] \leq N^{\epsilon}$ we found that 
\begin{equation}\label{eq-512}
\Sigma_{m_1,m_2,p} \ll   
\frac{NC(k_1,k_2)}{k_1[r_1^2 d_1 , r_2^2 d_2]\log^2 N} 
\prod_{\substack{ q\mid r_1r_2 \\ q > 3} }\left( 1- \frac{3}{q}\right)^{-1}.
\end{equation}
\end{case}
\begin{case}
For the case $N^{\epsilon} < [r_1^2 d_1 , r_2^2 d_2] < N^{2 - \epsilon}$ we use the trivial bound
\begin{equation*}
\Sigma_{m_1,m_2,p} \leq \left|\mathcal{A}\right| \leq \frac{N}{k_1 [T_1,T_2]} + 1.
\end{equation*}
From the definition of $T_1,T_2$ we see that $[T_1,T_2]\gg [r_1^2 d_1 , r_2^2 d_2]$ (with implied constants depending on $\Delta,D$).
And so, we get the bound 
\begin{equation}\label{eq-513}
\Sigma_{m_1,m_2,p} \leq \left|\mathcal{A}\right| \ll \frac{N}{k_1[r_1^2 d_1 , r_2^2 d_2]} + 1
\ll  \frac{N}{[r_1^2 d_1 , r_2^2 d_2]} + 1.
\end{equation}
\end{case}

\begin{case}
We now consider the case $[r_1^2 d_1 , r_2^2 d_2] \geq N^{2-\epsilon}$.
In this case, instead of bounding $\Sigma_{m_1,m_2,p}$ individually, we will give a bound for 
\begin{equation}\label{eq-555}
\mathop{\sum_{\substack{r_1,r_2\in\mathbb{N}\\ (r_1r_2,\Delta) = 1}} \;
\sum_{\substack{d_1,d_2 \\ q\mid d_1 d_2 \Rightarrow q \mid \Delta }}}_{[r_1^2d_1,r_2^2d_2] \geq N^{2-\epsilon}}\Sigma_{m_1,m_2,p}\;.
\end{equation}

Since $[r_1^2 d_1 , r_2^2 d_2] \geq N^{2-\epsilon}$ we must have that
$$
\text{max}(r_1^2 d_1,r_2^2 d_2)\geq  N^{1 - \epsilon/2}.
$$
Assume without loss of generality that $r_1^2 d_1 \geq N^{1 - \epsilon/2}$.
It follows that $m_1$ must satisfy $m_1\leq N^{\epsilon/2}$ (since otherwise $r_1^2 d_1 m_1 > N$ which is impossible).
For each such $m_1$ we want to count the number of $r_1,r_2,d_1,d_2,m_2,p$ which are counted in \eqref{eq-555}.
Note that each valid choice of $m_1,r_1,d_1$ determines $p,m_2,d_2,r_2$ uniquely.
And so, we get that
\begin{multline}\label{eq-556}
\mathop{\sum_{\substack{r_1,r_2\in\mathbb{N}\\ (r_1r_2,\Delta) = 1}} \sum_{\substack{d_1,d_2 \\ q\mid d_1 d_2 \Rightarrow q \mid \Delta }}}_{[r_1^2d_1,r_2^2d_2] \geq N^{2-\epsilon}}\Sigma_{m_1,m_2,p} \ll
\sum_{m_1 \leq N^{\epsilon/2}}
\sum_{r_1 \leq \sqrt{N}}
\sum_{\substack{d_1\leq N \\ q\mid d_1 \Rightarrow q \mid \Delta}} 1
\ll \\
N^{\epsilon/2}\sqrt{N} \prod_{q\mid \Delta}\log_q N
\ll N^{(1 + \epsilon)/2} \log^{\omega(\Delta)} N
\end{multline}
where $\omega(\Delta)$ is the number of distinct prime divisors of $\Delta$.
In particular, we have that 
\begin{equation}\label{eq-557}
\mathop{\sum_{\substack{r_1,r_2\in\mathbb{N}\\ (r_1r_2,\Delta) = 1}} \sum_{\substack{d_1,d_2 \\ q\mid d_1 d_2 \Rightarrow q \mid \Delta }}}_{[r_1^2d_1,r_2^2d_2] \geq N^{2-\epsilon}}\Sigma_{m_1,m_2,p} \ll N^{1 - \delta}
\end{equation}
for some $\delta > 0$.

\end{case}

We now split \eqref{eq-58} into three sums according to the cases considered above.
We get
\begin{multline*}
S_{k_1,k_2}(N) \ll \\
\mathop{\sum_{\substack{r_1,r_2 \\ (r_1r_2 ,\Delta) = 1}}\sum_{\substack{d_1,d_2 \\ q\mid d_1d_2 \Rightarrow q\mid \Delta}}}_{[r_1^2 d_1, r_2^2 d_2]\leq N^{\epsilon}} \Sigma_{m_1,m_2,p}
+ 
\mathop{\sum_{\substack{r_1,r_2 \\ (r_1r_2 ,\Delta) = 1}}\sum_{\substack{d_1,d_2 \\ q\mid d_1d_2 \Rightarrow q\mid \Delta}}}_{N^{\epsilon} < [r_1^2 d_1, r_2^2 d_2] < N^{2 - \epsilon}}\Sigma_{m_1,m_2,p}
+ \\
\mathop{\sum_{\substack{r_1,r_2 \\ (r_1r_2 ,\Delta) = 1}}\sum_{\substack{d_1,d_2 \\ q\mid d_1d_2 \Rightarrow q\mid \Delta}}}_{[r_1^2 d_1, r_2^2 d_2] \geq  N^{2 - \epsilon}} \Sigma_{m_1,m_2,p}.
\end{multline*}
Applying the bounds \eqref{eq-512}, \eqref{eq-513}, we get
\begin{multline}\label{eq-514}
S_{k_1,k_2}(N) \ll
\frac{N C(k_1,k_2)}{k_1\log^2 N} 
\left(\mathop{\sum_{\substack{r_1,r_2 \\ (r_1r_2 ,\Delta) = 1}}\sum_{\substack{d_1,d_2 \\ q\mid d_1d_2 \Rightarrow q\mid \Delta}}}_{[r_1^2 d_1, r_2^2 d_2]\leq N^{\epsilon}}\frac{\prod_{\substack{ q\mid r_1r_2 \\ q > 3}}\left(1 - \frac{3}{q}\right)^{-1}}{[r_1^2 d_1, r_2^2 d_2]}
\right)
+ \\
N 
\left(\mathop{\sum_{\substack{r_1,r_2 \\ (r_1r_2 ,\Delta) = 1}}\sum_{\substack{d_1,d_2 \\ q\mid d_1d_2 \Rightarrow q\mid \Delta}}}_{N^{\epsilon} < [r_1^2 d_1, r_2^2 d_2] < N^{2 - \epsilon}} \frac{1}{[r_1^2 d_1, r_2^2 d_2]}
\right)
+
\left(\mathop{\sum_{\substack{r_1,r_2 \\ (r_1r_2 ,\Delta) = 1}}\sum_{\substack{d_1,d_2 \\ q\mid d_1d_2 \Rightarrow q\mid \Delta}}}_{N^{\epsilon} < [r_1^2 d_1, r_2^2 d_2] < N^{2 - \epsilon}} 1
\right)
 +\\ 
\left(
\mathop{\sum_{\substack{r_1,r_2 \\ (r_1r_2 ,\Delta) = 1}}\sum_{\substack{d_1,d_2 \\ q\mid d_1d_2 \Rightarrow q\mid \Delta}}}_{[r_1^2 d_1, r_2^2 d_2] \geq  N^{2 - \epsilon}} \Sigma_{m_1,m_2, p}\right)
= \frac{N C(k_1,k_2)}{k_1\log^2 N} \Sigma_1 + N \Sigma_2 + \Sigma_3 + \Sigma_4.
\end{multline}

From \autoref{lem-a6} we get that $\Sigma_1 \ll 1$.
From \autoref{lem-a4} we get that that $\Sigma_2 \ll \frac{\log^2 N}{N^{\epsilon/2}}$.
From \autoref{lem-a5} we get that $\Sigma_3 \ll N^{1 - \epsilon/2}\log^2 N$.
Lastly, from \eqref{eq-557} we get that $\Sigma_4 \ll N^{1 - \delta}$.
Plugging these bounds into \eqref{eq-514} we get
\begin{equation}\label{eq-515}
S_{k_1,k_2}(N) \ll \frac{N C(k_1,k_2)}{k_1 \log^2 N}
\end{equation}
as required.

\end{proof}

\section{Further thoughts and conjectures}\label{sec-6}
\subsection{The asymptotics of $\psi_F$}
We proved that, under the conditions of \autoref{main-thm3}, the order of magnitude of $\psi_F(N)$ is $\frac{N \log \log N}{\sqrt{\log N}}$.
However, proving precise asymptotics 
$$
\psi_F(N) \sim c_F \frac{N \log \log N}{\sqrt{\log N}}
$$
is beyond our methods.

In \cite{Motohashi-1} Motohashi put forward a conjecture regarding the asymptotics for the number of primes up to $N$ which are of the form $x^2 + y^2 + 1$.
This conjecture was corrected by Iwaniec in \cite{Iwaniec-2}, to give
\begin{conjecture}
\begin{multline*}
C_1 := \lim_{N\rightarrow \infty}
\frac{\# \left\lbrace p \leq N \; \middle| \; p = x^2 + y^2 + 1\right\rbrace}{N \log^{3/2} N }= \\
\frac{1}{\sqrt{2}}\prod_{p\equiv 3 \bmod{4}}
\left(1 - \frac{1}{p^2}\right)^{-1/2}
\left(1 - \frac{1}{p(p-1)}\right)
\approx 0.610534\dots \;.
\end{multline*}
\end{conjecture}
We present some numerics to corroborate this conjecture.
For this, we compare the ratio
$$
R(N) = \left.\sum_{\substack{p \leq N \\ p = x^2 + y^2 + 1}} \log p \, \middle/ \sum_{\substack{n \leq N \\ n = x^2 + y^2 + 1}} 1  \right.$$
to $C_1 / L \approx 0.79889\dots $ where 
$$
L = \frac{1}{\sqrt{2}}\prod_{p\equiv 3 \bmod{4}}
\left(1 - \frac{1}{p^2}\right)^{-1/2} \approx 0.76422\dots$$
is the Landau-Ramanujan constant.
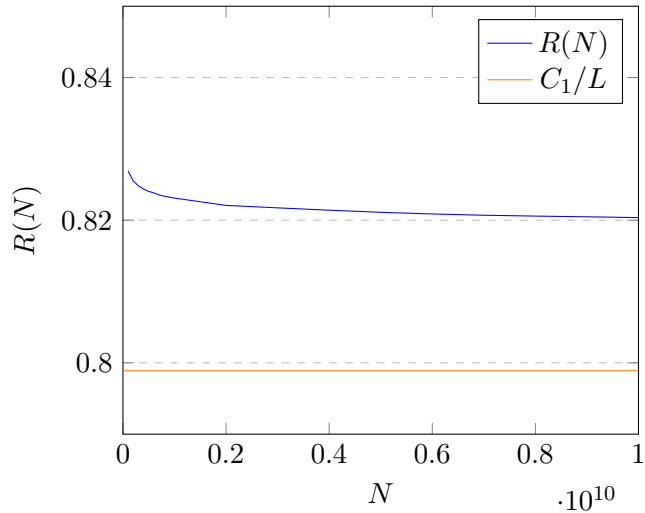
\begin{figure}[ht]
\begin{center}
\begin{tikzpicture}
\begin{axis}[
    xlabel={$N$},
    ylabel={$R(N)$},
    xmin=0, xmax=10000000000,
    ymin=0.79, ymax=0.85,
    legend pos=north east,
    ymajorgrids=true,
    grid style=dashed,
]
\addplot[
    color=blue,
    ]
    coordinates {
    (100000000, 0.8268908019101114)
 (200000000, 0.8255092426294897)
 (300000000, 0.8248157127445063)
 (400000000, 0.8243580149991577)
 (500000000, 0.8240407276868972)
 (600000000, 0.8238247661809822)
 (700000000, 0.8235507305730502)
 (800000000, 0.8233758840110366)
 (900000000, 0.8232186845862611)
 (1000000000, 0.8230996308570618)
 (2000000000, 0.8220738017016627)
 (3000000000, 0.8217295070029147)
 (4000000000, 0.8214044778532162)
 (5000000000, 0.8211043169753866)
 (6000000000, 0.8208739173315096)
 (7000000000, 0.8206940976651872)
 (8000000000, 0.8205762221458976)
 (9000000000, 0.8204763457778087)
 (10000000000, 0.8203651242177011)
    };
    \addlegendentry{$R(N)$}
\addplot[
    color=orange,
    ]
    coordinates {
    (0, 0.7988946498397742)
    (10000000000, 0.7988946498397742)
    };
    \addlegendentry{$C_1 / L$}
\end{axis}
\end{tikzpicture}
\end{center}
\caption{ Plot of $R(N)$, the ratio of the prime counting function $\sum_{p=x^2+y^2+1\leq N} \log p$ to the number of integers $n\leq N$ of the form $n=x^2+y^2+1$.}
\end{figure}\label{fig1}
\autoref{fig1} shows that there is an error of approximately $2.7\%$ for $N = 10^{10}$.

Following the same heuristics, we conjecture that
$$
C_k := \lim_{N\rightarrow \infty}
\frac{\# \left\lbrace p \leq N/k \; \middle| \; kp = x^2 + y^2 + 1\right\rbrace}{N \log^{3/2} N } 
$$
exists and satisfies
\begin{equation*}
C_k =
\begin{cases}
 C_1 \frac{(k,2)}{k}
\prod_{\substack{p\mid k \\ p\equiv 3 \bmod{4}}}\frac{p^2 - 1}{p^2 - p - 1} & 4\nmid k \\ 
0 & 4\mid k
\end{cases}.
\end{equation*}

With these assumptions, we can give a conjecture for the constant arising in the asymptotics of $\psi_{x^2 + y^2 + 1}(N)$:
\begin{conjecture}
For $F = x^2 + y^2 + 1$, 
$$
\psi_F(N) \sim c_F \frac{N \log \log N}{\sqrt{\log N}}
$$
with
\begin{multline*}
c_F = \lim_{N\rightarrow\infty} \frac{1}{\log N}\sum_{k\leq N}C_k
= \\
\frac{1}{\sqrt{2}}\prod_{p\equiv 3 \bmod{4}}
\left(1 - \frac{1}{p^2}\right)^{-1/2}
\left(1 - \frac{1}{p(p-1)}\right)
\left(1 + \frac{1}{p^2 - p - 1}\right) = \\
\frac{1}{\sqrt{2}}\prod_{p\equiv 3 \bmod{4}}
\left(1 - \frac{1}{p^2}\right)^{-1/2} = L \approx 0.76422 \dots \;.
\end{multline*}
\end{conjecture}

\subsection{LCM for $x,y$ in a box}
We have chosen to consider 
$$
\psi_F(N) = \log\left(
\mathop{\text{LCM}}_{0<F(x,y)\leq N}
\left\lbrace F(x,y)\right\rbrace
\right)
$$
as a generalization of the LCM problem in the single variable case.
Another possible generalization is given by
$$
L_F(N) = \log\left(
\mathop{\text{LCM}}_{|x|,|y|\leq \sqrt{N} }
\left\lbrace F(x,y)\right\rbrace
\right).
$$
In the case of \autoref{main-thm3}, $L_F(N)$ has the same order of magnitude as $\psi_F(N)$.
We do not give a detailed proof of this statement, but only highlight the main idea for the equivalence.

\autoref{lem-ap1} allows us to reduce the problem of computing $\psi_F(N)$ for a general polynomial $F(x,y)$ to computing $\psi_F(N)$ for $F$ a shifted quadratic form $\varphi(x,y) + C$.
So it is enough to show that $\Theta\left(L_F(N)\right) = \Theta\left(\psi_F(N)\right)$ for the case $F(x,y)=\varphi(x,y) +C $

If $\varphi$ is positive definite, then the equivalence is obvious.
In the indefinite case, it follows from reduction theory for quadratic forms that if $n$ is represented by $\varphi$ then it is possible to find a representation $n=\varphi(x,y)$ with $|x|,|y|\ll \sqrt{n}$.
From this, one can once more show that $\Theta\left(L_F(N)\right) = \Theta\left(\psi_F(N)\right)$.



\end{document}